\documentclass[12pt,reqno]{amsart}

\usepackage{amscd,amssymb}

\newtheorem{theorem}{Theorem}[section]
\newtheorem{proposition}[theorem]{Proposition}
\newtheorem{lemma}[theorem]{Lemma}
\newtheorem{corollary}[theorem]{Corollary}

\theoremstyle{definition}
\newtheorem{definition}[theorem]{Definition}

\theoremstyle{remark}

\numberwithin{equation}{section}

\newcommand{\Z}{{\mathbf Z}}
\newcommand{\Q}{{\mathbf Q}}
\newcommand{\R}{{\mathbf R}}
\newcommand{\C}{{\mathbf C}}
\renewcommand{\H}{{\mathbb H}}

\newcommand{\OO}{{\mathcal {O}}}
\newcommand{\RR}{{\mathcal {R}}}

\newcommand{\M}{{\mathcal {M}}}
\newcommand{\T}{{\mathcal {T}}}
\newcommand{\m}{{\mathfrak {m}}}
\renewcommand{\L}{{\mathcal {L}}}

\DeclareMathOperator{\id}{{id}}

\DeclareMathOperator{\Hom}{{\operatorname{Hom}}}
\DeclareMathOperator{\im}{\operatorname{im}}

\DeclareMathOperator\X{{\mathcal X}}
   
\DeclareMathOperator\Y{{\mathcal Y}}  

\DeclareMathOperator\e{{\mathfrak e}} 
\renewcommand\c{{\mathfrak c}}

\DeclareMathOperator\tr{\operatorname{tr}} 

\DeclareMathOperator\ns{\operatorname{ns}}     
\DeclareMathOperator\ob{\operatorname{Ob}}    
\DeclareMathOperator\cl{\operatorname{cl}} 
\renewcommand\H{\mathcal H}
\DeclareMathOperator\E{\mathcal E}
\def\<{\langle}
\def\>{\rangle}

\DeclareMathOperator\CC{\mathcal C}
\DeclareMathOperator\eca{\mathcal E({\mathcal C})}

\renewcommand\L{\mathcal L}
\DeclareMathOperator\lam{\lambda}

\DeclareMathOperator\I{\mathcal I}

\DeclareMathOperator\tdim {{\mathfrak {tordim}}}

\DeclareMathOperator\lo{{{\text{Lim}_{\omega}}}}

\DeclareMathOperator\ex{{\mathcal Exc}}
\DeclareMathOperator\ts{\tilde\otimes_\Lambda}
\DeclareMathOperator\tsig{{\mathfrak {torsign}}}
\DeclareMathOperator \df{\mathfrak D}
\DeclareMathOperator\s{\mathfrak s}

\DeclareMathOperator\coker{\operatorname{coker}}

\begin{document}

\title[Novikov - Shubin signatures, II]
{Novikov - Shubin signatures, II}

\author[M.~Farber]{M.~Farber}
\address{Department of Mathematics, Tel Aviv University, Tel Aviv, 69978, Israel}
\email{farber@math.tau.ac.il}
\date{August 1999}

\thanks{The research was partially supported by a grant from the 
Israel Academy of Sciences and Humanities and by
the Herman Minkowski Center for Geometry}

\begin{abstract}
This paper continues \cite{Fa3}. Here we construct a linking form on the torsion part 
of middle dimensional extended $L^2$ homology and cohomology of odd-dimensional manifolds.
We give a geometric necessary  condition when this linking form is hyperbolic.
We compute this linking form in case, when the manifold bounds. 
We introduce and study new numerical invariants of the linking form: 
the Novikov - Shubin signature and the torsion signature; we compute these invariants
explicitly for manifolds with $\pi_1=\mathbf Z$ in terms of the Blanchfield form. We develop
a notion of excess for extensions of torsion modules and show how this concept 
can be used to guarantee vanishing of the torsion signature.
\end{abstract}

\maketitle

\section{\bf  Introduction}\label{intr}

In the classical topology of manifolds it is well known that the torsion subgroup $T$
of one-dimensional
integral homology of a 3-dimensional oriented manifold supports a symmetric nondegenerate
linking form $$T\times T\to \Q/\Z.$$ 
Equivariant analogs of the linking form (for example, the Blanchfield form)
play an important role in the knot theory.

The purpose of the present paper is to construct a linking form of a new kind, 
which lives on the torsion part of the extended $L^2$-homology. 
Namely, let $M$ be a closed oriented
manifold of dimension $2q+1$. For any flat Hermitian bundle $\E$ over $M$ with fiber and monodromy in a finite von Neumann category $\CC$, we defined in \cite{Fa}, \cite{Fa1}, \cite{Fa2},
\cite{Fa4} extended $L^2$ homology $\H_i(M;\E)$, which is an object of an abelian category $\eca$ determined by $\CC$. 
The homology $\H_i(M;\E)$ splits canonically as a direct sum of its projective and torsion parts.
We show in this paper that the torsion part of the middle dimensional homology $T(\H_q(M;\E))$
supports a canonical linking form which carries a topological information.

Classification of the 
torsion Hermitian forms of this sort was performed in the previous paper \cite{Fa3}. 

In the present paper we study
numerical invariants of such linking forms, the Novikov - Shubin signature and the torsion signature.
These invariants use an additional structure, the trace on category $\CC$, and for the torsion signature
we need this trace to be not normal, i.e. Dixmier type. The Novikov - Shubin signature 
is based on the Novikov - Shubin invariant; 
the torsion signature is based on the concept of the tosion dimension
which was introduced by the author in \cite{Fa4}. 

We compute the Novikov - Shubin signature and the torsion signature in the
(commutative) case of
manifolds with fundamental group $\Z$; we show that in this case they can be expressed completely
in terms of the Blanchfield form. 

The importance of the torsion signature and the Novikov - Shubin
signature follows from the fact that these invariants are defined in a 
general non-commutative situation.  We should also emphasize existence of a 
very large variety of von Neumann 
flat Hermitian bundles over a given manifold; we refer to \cite{Fa4}, where such bundles are 
associated with growth processes. 

I hope that the algebraic invariants of linking forms, in the spirit of the ones
constructed in the present paper, will provide
useful topological invariants of manifolds. For instance, 
these invariants can be applied to study knots and links. 
I plan to describe applications to the knot concordance problem in a separate article.

\section{ \bf Push forward Hermitian forms}\label{sec1}

In this section we describe a technique, which allows to apply functors to Hermitian forms.
Our major motivation is as follows. Starting from a closed PL manifold, we have 
the Poincar\'e duality, which may be viewed as a Hermitian form in the homotopy 
category of chain complexes (cf. section \ref{sec1}); 
the push forward forms of the Poincar\'e duality will be 
different intersection and linking forms. One of these (the 
linking form on the torsion part of the extended
$L^2$ cohomology) is the main subject of this paper.

\subsection{}\label{sec1.1}
We will use the formalism of Hermitian forms in categories with duality as described in \cite{Fa3}, \S 1.

Suppose that $\CC$ and $\CC^\prime$ are two categories with dualities
$(D, s)$ and $(D^\prime, s^\prime)$, correspondingly.
Let 
$F:\CC \to \CC^\prime$ be a covariant functor. We will suppose that $F$ has the following property.
For any symmetry $\epsilon\in G(\CC)$ (cf. \cite{Fa3}, \S 1.2) there
exists a unique symmetry $\epsilon'\in G(\CC^\prime)$, such that for any $M\in\ob(\CC)$
holds $F(\epsilon_M) = \epsilon'_{F(M)}.$ In this situation $F$ determines a homomorphism
$$F_\ast: G(\CC)\to G(\CC^\prime), \quad F_\ast(\epsilon) =\epsilon'$$
between the groups of symmetries.

The functor $F$ determines also two contravariant functors $FD, \, \, \, D^\prime F: \CC\to \CC^\prime$. 

Let $\gamma: FD\to D^\prime F$ be a natural transformation.

Given a form
$\phi:M\to D(M)$ in category with duality $\CC$ (cf. \cite{Fa3}, \S 1.3), 
consider the following form in $\CC^\prime$:
\begin{eqnarray}
\psi: F(M) \to D^\prime F(M),\quad\text{where}\quad
\psi = \gamma(M)\circ F(\phi).\label{(1-1)}
\end{eqnarray} 
We will say that form $\psi$ is {\it derived} from $\phi$ by means of the pair $(F, \gamma)$, or 
that it is {\it a push forward of} $\phi$. 

If $\gamma$ is a natural isomorphism of functors,
then the push forward form $\psi$ of any non-degenerate form $\phi$
is non-degenerate.

We wish to understand under which conditions the push forward form $\psi$ is Hermitian.

\begin{definition} We will say that the pair $(F,\gamma)$ 
is $\eta$-{\it Hermitian}, where $\eta\in G(\CC^\prime)$, 
if the following diagram of natural transformations
\begin{eqnarray}
\CD
F@>{Fs}>>FDD\\
@V{\eta s^\prime F}VV  @VV{\gamma D}V\\
D^\prime D^\prime F@>>{D^\prime\gamma}>D^\prime FD
\endCD\label{(1-2)}
\end{eqnarray}
commutes. 
\end{definition}
Here $Fs$ denotes the natural transformation between the
functors $F$ and $FDD$ which to any object $M\in\ob(\CC^\prime)$ associates 
the morphism $F(s(M)): F(M)\to $ $FDD(M)$. Similarly, 
the natural transformation $\gamma D: FDD \to D^\prime FD$
associates to any $M\in\ob(\CC^\prime)$ the morphism
$\gamma(D(M)): FDD(M) \to D^\prime FD(M)$. This explains our
notations.

\begin{lemma}\label{l1.3}
Suppose that $(F, \gamma)$ 
is $\eta$-Hermitian. If the initial form $\phi: M\to D(M)$ is $\epsilon$-Hermitian, 
then the push forward form $\psi: F(M)\to D^\prime F(M)$ in $\CC^\prime$  is $\epsilon'$-Hermitian, 
where $\epsilon' = \eta^{-1}F_\ast(\epsilon)\in G(\CC^\prime)$.\end{lemma}
\begin{proof} Consider the following diagram
\begin{eqnarray}
\CD
F(M)@>{F(s(M))}>>FD D(M)@>{FD(\phi)}>>FD(M)\\
@V{\eta s^\prime(F(M))}VV @V{\gamma(D(M))}VV @VV{\gamma(M)}V\\
D^\prime D^\prime F(M)@>>{D^\prime(\gamma(M))}>D^\prime FD(M)@>>{D^\prime F(\phi)}>D^\prime F(M).\label{(1-3)}
\endCD
\end{eqnarray}
The square on the right is commutative because $\gamma$ is a natural 
transformation; the square on the left is commutative since (\ref{(1-2)}) 
holds. One of the compositions from the left upper corner of (\ref{(1-3)}) to 
the right lower corner is
$$ \eta\circ D^\prime F(\phi)\circ D^\prime(\gamma(M))\circ
s^\prime(F(M))= \eta\psi^\dagger.$$
The other composition equals
$$\gamma(M)\circ F(\phi^\dagger) = \gamma(M)\circ F(\epsilon\phi) =
F_\ast(\epsilon)\gamma(M)\circ F(\phi) = F_\ast(\epsilon)\psi.$$
\end{proof}

\section{\bf  Poinca\'e duality as a Hermitian form}\label{sec2}

\subsection{Duality in the category of chain complexes} 
Here we will recall the construction of a duality in the category of chain complexes,
which geometrically corresponds to the Poincar\'e duality for manifolds. 

Let $\CC$ be an additive category
with duality $(D,s)$.  Fix a non-negative
integer $n$ and consider the category $\CC_n(\CC)$, which has as its objects chain complexes 
of fixed length $n+1$ 
$$C\ =\ (0\to C_n\to \dots C_{i+1}\stackrel{d}\to C_i\stackrel{d}\to C_{i-1}\to \dots\to C_0\to 0)
$$
over $\CC$. Morphisms of this category are homotopy classes of 
$\CC$-chain morphisms.

We will describe a duality in $\CC_n(\CC)$.
For any chain complex $C\in \ob(\CC_n(\CC))$ denote by $\df (C)$ the following 
chain complex
$$\df (C)\ =\ (0\to D_n\stackrel{(-1)^n d^\ast}\to D_{n-1}\dots \to 
D_i\stackrel{(-1)^{i}d^\ast}\to
D_{i-1}\dots \to D_0\to 0),$$
where $D_i = D(C_{n-i})= C_{n-i}^\ast$ is the dual of 
$C_{n-i}$ with respect to
the duality $(D,s)$ in $\CC$, and $d^\ast:D_i\to D_{i-1}$ is 
$d^\ast=D(d):D(C_{n-i})\to D(C_{n-i+1})$, the dual
of the boundary homomorphisms
$d:C_{n-i+1}\to C_{n-i}$ (with respect to $(D,s)$). 

If $f:C\to C'$ is a chain map between chain complexes in $\CC_n(\CC)$, which is
given by $\CC$-morphisms $f_i:C_i\to C'_i$, for $i=0,1,\dots, n$,
then the dual morphism $\df (f): \df (C')\to \df (C)$ is given by the sequence
$D(f_i):D(C'_i)\to D(C_i)$. Note also that the homotopy class of $\df (f)$
is determined by the homotopy class of $f$.
Thus we obtain a contravariant functor $\df :\CC_n(\CC) \to\CC_n(\CC)$.

The natural isomorphism 
$\s:\id_{\CC_n(\CC)} \to \df \circ \df ,$
which we have to specify in order to define a duality in the category $\CC_n(\CC)$ 
(cf. \cite{Fa3}, \S 1.1)
assigns to any chain
complex $C$ the chain map $\s(C): C\to \df \df (C)$, which is given by the
isomorphisms $\epsilon_i s_i: C_i\to DD(C_i)$, where $i=0,1,2,\dots, n$,
where $s_i=s(C_i)$, and 
$\epsilon_i = (-1)^{i(n+1)}.$

The above construction applies to the category $\Lambda$-mod 
of finitely generated projective left $\Lambda$-modules,
where $\Lambda$ is the ring with involution, with its canonical duality, cf. \cite{Fa3}, \S 1.4.
The corresponding homotopy
category of chain complexes $\CC_n(\Lambda\text{-mod})$ 
will be briefly denoted $\CC_n(\Lambda)$.

\subsection{Poincar\'e duality}
Here we recall the classical                                 
construction, which associates to any $n$-dimensional manifold a 
Hermitian form in the homotopy category of chain complexes
$\CC_n(\Lambda)$.

Let $K$ be a closed piecewise linear manifold of dimension $n$. We will
denote by $\pi$ its fundamental group $\pi=\pi_1(K)$. Let $w:\pi\to \{1,-1\}$
be the orientation character of $K$ (the first Stiefel - Whitney class).
It determines an involution on the group ring $\Lambda=\Q[\pi]$ given by
$g\mapsto \overline g = w(g)g^{-1}$ for $g\in \pi$. 

Fix two mutually dual triangulations on $K$ and denote by $C(\tilde K)$ and $C'(\tilde K)$
the simplicial chain complexes of the universal covering $\tilde K$ with
respect to these triangulations. $C(\tilde K)$ and $C'(\tilde K)$ consist of free $\Lambda$-modules and have length 
$n+1$ and so they can
be viewed as objects of category $\CC_n(\Lambda)$. Any choice of orientation of the universal covering $\tilde K$ determines a non-degenerate 
{\it intersection pairing} (cf. \cite{B}, \cite{M2})
\begin{eqnarray}
[\ ,\ ]: C(\tilde K) \times C'(\tilde K) \to \Lambda,\label{(2-1)}
\end{eqnarray}
which is $\Lambda$-linear with respect to the first variable and anti-linear
with respect to the second variable. It has the property
\begin{eqnarray}
[dx,y] = (-1)^{p}[x,dy]
\label{(2-2)}
\end{eqnarray}
for $x\in C(\tilde K)$ and $y\in C'(\tilde K)$, where $p=|x|$ denotes the dimension of $x$.
Note that $[x,y]$ is 0 unless $|x| + |y| = n$. Recall that $[\, ,\, ]$ is defined by the formula
$$[x,y] = \sum_{g\in \pi} g^{-1}\langle (gx),y\rangle,$$
where $\langle x,y\rangle$ denotes the usual geometric intersection number of chains $x$ and $y$
with respect to a given orientation of $\tilde K$ (adding signs at each intersection point).

Since $C(\tilde K)$ and $C'(\tilde K)$ are constructed with respect to two different
triangulation of the same space, there is a chain homotopy equivalence
$h:C(\tilde K)\to C'(\tilde K)$, which is determined uniquely up to homotopy. In fact (since $\frac{1}{2}\in \Lambda)$) $h$
can be chosen so that
\begin{eqnarray}
[x,h(y)] = (-1)^{pq}\overline{[y,h(x)]},\quad\text{where}\quad 
p=|x|, \, \, q=|y|,\label{(2-3)}
\end{eqnarray}
for all $x, y \in C(\tilde K)$. 

The map
$x\mapsto (y\mapsto [x,h(y)]\in \Lambda)$
determines a chain map
\begin{eqnarray}
\phi: C(\tilde K) \to \df(C(\tilde K)),\quad \phi(x)(y)=[x,h(y)]\label{(2-4)}
\end{eqnarray}
(this statement is equivalent to (\ref{(2-2)}), which we will view as a 
non-degenerate form in 
category $\CC(\Lambda)$ supplied with duality $(\df,\s)$. Note that the
canonical isomorphism $\s$ in $\CC_n(\Lambda)$ is constructed as explained in section 2.1.
An easy check using the definitions shows that
\begin{eqnarray*}\phi^\dagger = \phi,\label{(2-5)}\end{eqnarray*}
i.e. $\phi$ is Hermitian.
\begin{corollary} \label{cor2.3}
Any closed piecewise linear manifold $K$ of dimension $n$ with oriented universal cover $\tilde K$
determines canonically (via (\ref{(2-4)}) a nondegenerate Hermitian form in 
category $\CC_n(\Lambda)$ 
with respect to duality $(\df,\s)$. 
\end{corollary}

\section{\bf Linking form on extended $L^2$-homology}\label{sec3}

In this section we introduce the linking form on the extended $L^2$-homology and
cohomology. 

\subsection{Categories $\Lambda$-mod-$\CC$ and $\CC$-mod-$\Lambda$} 
Let $\CC$ be a Hilbertian von Neumann
category (cf. \cite{Fa3}, \S 2.2). If $\Lambda$ is a ring, we may
consider the category of left $\Lambda$-modules in $\CC$, cf. \cite{Fa2}, \S 6. An object of $\Lambda$-mod-$\CC$ is defined as $\X\in\ob(\CC)$ with
a given ring homomorphism $\Lambda\to \hom_{\CC}(\X,\X)$. We may think of
$\Lambda$ as acting on $\X$ from the left. Morphisms of 
$\Lambda$-mod-$\CC$ are the morphisms of $\CC$ which commute with the 
action of $\Lambda$. It is clear that $\Lambda$-mod-$\CC$ is an additive
category. 

Similarly one may consider the category $\CC$-mod-$\Lambda$ of {\it right} 
$\Lambda$-modules in $\CC$. A right $\Lambda$-module in $\CC$ is an object
$\X$ of $\CC$ with a ring homomorphism $\Lambda^{op}\to \hom_{\CC}(\X,\X)$,
where $\Lambda^{op}$ is the opposite ring of $\Lambda$. 

If $\M$ is a left $\Lambda$-module in $\CC$ then the dual $\M^\ast$ in $\CC$
(cf. \cite{Fa3}, \S 2) is
naturally defined as a right $\Lambda$-module:
\begin{eqnarray*}
(\phi\cdot\lambda)(m) = \phi(\lambda m)\quad\text{for}\quad
\lambda\in \Lambda, m\in \M.
\end{eqnarray*}
Conversely, the dual in $\CC$ of a right $\Lambda$-module is a left
$\Lambda$-module.

Suppose that we are given an involution on $\Lambda$. Then we may canonically
construct a duality $(D_\Lambda,s_\Lambda)$ in $\Lambda$-mod-$\CC$. 
Namely, given a $\Lambda$-module
$\M$ in $\CC$, we will define $D_\Lambda(\M)$ to be the dual of $\M$ in $\CC$
(cf. 2.2 in \cite{Fa3}, i.e. $D_\Lambda(\M)$ is the set of all anti-linear functionals on $\M$)
with the following action of $\Lambda$. If $\phi\in D_\Lambda(\M)$ and 
$\lambda\in \Lambda$ we set 
\begin{eqnarray*}
(\lambda\cdot\phi)(m) = \phi(\overline \lambda\cdot m),
\quad\text{for}\quad
m\in \M.\end{eqnarray*} 
One checks that the canonical isomorphism 
$s_\Lambda=s: \M \to D_\Lambda D_\Lambda(\M)$ in $\CC$ 
is now an isomorphism in $\Lambda$-mod-$\CC$, i.e. it commutes with the
$\Lambda$-action. 
Let us emphasize that the duality in $\Lambda$-mod-$\CC$ 
depends on the involution of $\Lambda$ in an essential way.

Similarly, one introduces a duality in 
$\CC$-mod-$\Lambda$ using an involution of the ring $\Lambda$.

\subsection{} Let $K$ be a closed PL manifold of dimension $n$. Let $C(\tilde K)$ denote the 
simplicial chain
complex of the universal covering $\tilde K$. From \S 2 we know that $C(\tilde K)$ carries a nondegenerate
Hermitian form $\phi: C(\tilde K)\to \df(C(\tilde K))$ in $\CC_n(\Lambda)$, 
whose sign depends on the choice of orientation of $\tilde K$. Here $\Lambda=\Q[\pi]$, where
$\pi=\pi_1(K)$,
considered together with the involution determined by the first Stiefel-Whitney class $w:\pi\to \{1,-1\}$.

Let $\CC$ be a Hilbertian von Neumann category. 
Let $\M\in \ob(\CC\text{-mod-}\Lambda)$ be a right $\Lambda$-module in $\CC$, and let 
$\psi: \M\to \M^\ast$ be an Hermitian form in $\ob(\CC\text{-mod-}\Lambda)$.
$\M$ is simply a unitary representation of $\pi$ with respect to the scalar product determined by
$\psi$. More precisely, the scalar product $\langle x,y\rangle =\psi(x)(y)$ for $x, y\in \M$ and the
right action of $\pi$ on $\M$ satisfy
$$\langle xg,yg\rangle = w(g)\langle x,y\rangle, \quad g\in \pi,\quad \text{and}\quad
\langle x,y\rangle = \overline{\langle y,x\rangle}.$$

The following composite
\begin{eqnarray}
\sigma: \M \tilde\otimes_\pi C(\tilde K)  \stackrel{\psi\otimes\phi}\longrightarrow \M^\ast
\tilde \otimes_\pi \df(C(\tilde K)) \stackrel{\otimes}\longrightarrow 
\df(\M\tilde \otimes_\pi C(\tilde K))
\label{(3-3)}\end{eqnarray}
can be viewed as a Hermitian form in category $\CC_n(\CC)$. 
We will use the following Lemma and
the push forward construction from section \ref{sec1} 
in order to define derived Hermitian forms of (\ref{(3-3)}).

\begin{lemma}\label{l3.3} 
Let $\CC$ be a Hilbertian von Neumann category and let $\T$
denote the torsion subcategory of the extended category $\eca$. Let
 $F:\CC_{2q+1}(\CC) \to \T$ be the functor which assigns to a chain 
complex $C\in\ob(\CC_{2q+1}(\CC))$ the torsion part of its $q$-dimensional
extended homology
\begin{eqnarray*}F(C) = T(\H_{q}(C))\in\ob(\T).\end{eqnarray*}
Consider the duality $(\df ,\s)$ in $\CC_{2q+1}(\CC)$ (cf. 3.1)
and the canonical duality 
$(\e,s)$ in $\T$,
cf. \cite{Fa3}, \S 3.3. Then there exists a natural equivalence
\begin{eqnarray*}\gamma: F\df \to \e F,\end{eqnarray*}  
so that the pair $(F, \gamma)$ is $(-1)^{q+1}$-Hermitian (in the sense of \S 1).
\end{lemma}
Proof will be given later in subsection \ref{sec3.6}.
\begin{theorem}\label{thm3.4} 
Let $K$ be a closed piecewise linear manifold of odd dimension $2q+1$.
Let $\M\in \ob(\CC\text{-mod-}\Lambda)$ be a representation of $\pi=\pi_1(K)$, where 
$\CC$ is a Hilbertian von Neumann category.
Let $T_{q}$ denote the torsion part of the extended $L^2$-homology
$T_{q}=T(\H_{q}(K;\M))$. Then a choice of orientation of the universal covering 
$\tilde K$ and of a Hermitian form $\psi: \M\to \M^\ast$ in 
$\CC\text{-mod-}\Lambda$ determine canonically
a non-degenerate $(-1)^{q+1}$-Hermitian form
\begin{eqnarray}\L: T_{q}\to \e(T_{q}).\label{(3-6)}\end{eqnarray}
\end{theorem}

This Theorem follows immediately from Lemma \ref{l1.3}, Corollary \ref{cor2.3} and Lemma \ref{l3.3}.

Form $\L$ above will be called {\it the linking form.}

\subsection{Proof of Lemma \ref{l3.3}}\label{sec3.6} First, note that the functor 
$F$ applied to a 
chain complex $C$ in $\CC_{2q+1}(\CC)$ gives the following torsion object
\begin{eqnarray*}F(C) = (d: C_{q+1}/Z_{q+1} \to \overline B_{q}),\end{eqnarray*}
where $Z_{q+1}$ denotes the space of cycles in $C_{q+1}$ and $B_{q}$ 
denote
the subspace of boundaries in $C_{q}$; cf. \cite{Fa1}, \cite{Fa2}. 
Applying the duality functor
$\e$ (cf. \cite{Fa3}, \S 3.3) we obtain
\begin{eqnarray*}\e(F(C)) = (d^\ast: (\overline B_{q})^\ast\to 
(C_{q+1}/Z_{q+1})^\ast),\end{eqnarray*}
where the star denotes the duality functor in $\CC$, cf. \cite{Fa3}, \S 2.2.   

Now we want to compute $F(\df(C))$. The chain complex $\df (C)$ in 
dimensions $q+1$ and $q$ looks as follows
\begin{eqnarray*}\dots \to C_{q}^\ast \stackrel{(-1)^{q+1}d^\ast}\longrightarrow C_{q+1}^\ast\to\dots.
\end{eqnarray*} 
Clearly $\cl(\im(d^\ast)) = \{\phi\in C^\ast_{q+1};\phi|_{Z_{q+1}}=0\} =
(C_{q+1}/Z_{q+1})^\ast$
and similarly
$$\ker(d^\ast) = \{\phi\in C^\ast_{q};\phi|_{\overline B_{q}}=0\},
\quad \text{and thus}\quad C_{q}^\ast/\ker(d^\ast)\simeq 
(\overline B_{q})^\ast.$$
Hence we obtain that $F(\df (C))$ can be identified with
\begin{eqnarray*}F(\df (C)) = ((-1)^{q+1}d^\ast : (\overline B_{q})^\ast 
\to (C_{q+1}/Z_{q+1})^\ast).\end{eqnarray*}
Therefore, we may define the natural transformation 
$\gamma(C): F\df (C) \to \e F(C)$ as the morphism in the 
extended category $\eca$ represented by the following commutative diagram
\begin{eqnarray*}
\CD
((-1)^{q+1}d^\ast: & (\overline B_{q})^\ast @>>> (C_{q+1}/Z_{q+1})^\ast)\\
& @V{(-1)^{q+1}\id}VV                      @VV{\id}V\\
(d^\ast : & (\overline B_{q})^\ast @>>> (C_{q+1}/Z_{q+1})^\ast),
\endCD
\end{eqnarray*}
which is clearly an isomorphism in $\eca$. 

Now examining diagram (\ref{(1-2)}) 
we find that it is commutative with $\eta=(-1)^{q+1}$ and hence the pair
$(F, \gamma)$ is $(-1)^{q+1}$-Hermitian.
\qed

\subsection{Computing the linking form}\label{sec3.7}
Here we will show how 
one may practically compute the linking form (\ref{(3-6)}). 

We will assume that $K$ is a odd-dimensional manifold, $n=\dim K=2q+1$. 
$\pi$ will denote the fundamental group of $K$, $w: \pi\to \Z_2$ the
orientation character, and the group ring $\Lambda=\Q[\pi]$ will be 
considered with the involution determined by $w$.

Fix two mutually dual triangulations of $K$ and consider the chain
complexes of the universal cover $\tilde K$ which we will denote by
$C_\ast(\tilde K)$ and $C'_\ast(\tilde K)$. Now, we have the following
commutative diagram
\begin{equation}
\begin{array}{ccc}
C_{q+1}(\tilde K) & \stackrel d\to & C_q(\tilde K)\\
{h_{q+1}}\downarrow & & \downarrow {h_q} \\
C'_{q+1}(\tilde K) & \stackrel{d'}\to & C'_q(\tilde K)
\end{array}\label{(3-12)}
\end{equation}
where $h_\ast$ denotes a {\it symmetric} (i.e. satisfying (\ref{(2-3)})
lift to the universal cover of a simplicial 
approximation of the identity map $K\to K$. Also, we have the following
two non-degenerate pairings (cf. (\ref{(2-1)}))
\begin{eqnarray*}C_{q+1}(\tilde K) \times C'_q(\tilde K) \to \Lambda
\label{(3-13)}\end{eqnarray*}
and 
\begin{eqnarray*}C_{q}(\tilde K) \times C'_{q+1}(\tilde K) \to \Lambda
\label{(3-14)}\end{eqnarray*}   
which are linear with respect to the first variable and anti-linear with 
respect to the second variable. They determine isomorphisms
\begin{eqnarray*}C'_{q+1}(\tilde K) \stackrel{\simeq}\longrightarrow D(C_q(\tilde K)), \end{eqnarray*}   
\begin{eqnarray*}C'_{q}(\tilde K) \stackrel{\simeq}\longrightarrow D(C_{q+1}(\tilde K)),\end{eqnarray*}   
where $D$ denotes the duality functor in the category of projective 
$\Lambda$-modules (cf. \cite{Fa3}, \S 1.4). We have
\begin{eqnarray*}
d' = (-1)^{q+1} D(d)\quad\text{and, using (\ref{(2-3)}),}\quad h_{q+1} = (-1)^{q(q+1)}D(h_q) = 
D(h_q)
\end{eqnarray*}
and thus we may rewrite diagram (\ref{(3-12)}) as follows
\begin{eqnarray}
\CD
C_{q+1}(\tilde K) @>d>> C_q(\tilde K)\\
@V{D(h)}VV  @VV{h}V \\
D(C_{q}(\tilde K)) @>{(-1)^{q+1}D(d)}>> D(C_{q+1}(\tilde K)).
\endCD\label{(3-18)}
\end{eqnarray}
Now, if we are given a right $\Lambda$-module $\M$ in a Hilbertian category
$\CC$ with an Hermitian non-degenerate form 
$\psi: \M\to D_\Lambda(\M)$, applying the functor $\ts$ (cf. \cite{Fa2}, \S 6.4) to diagram
(\ref{(3-18)}, we obtain the 
following commutative diagram
\begin{eqnarray}
\CD
\M\ts C_{q+1}(\tilde K) @>d>> \M\ts C_q(\tilde K)\\ 
@V{(-1)^{q+1} D(h')}VV  @VV{h'}V \\    
(\M\ts C_{q}(\tilde K))^\ast @>d^\ast>> (\M\ts C_{q+1}(\tilde K))^\ast ,
\endCD
\label{(3-19)}
\end{eqnarray}
where $h'$ denotes the composition
$$\M\ts C_q(\tilde K) \stackrel{\psi\otimes h}\longrightarrow 
D_\Lambda(\M)\ts D(C_{q+1}(\tilde K)) \simeq (\M\ts C_{q+1}(\tilde K))^\ast.$$ 
The last diagram (\ref{(3-19)} clearly represents a $(-1)^{q+1}$-Hermitian
form on the torsion part $T(\coker d) = T(\H_q(K,\M))$ (compare \cite{Fa3}, \S 4), 
which coincides with
the linking form (\ref{(3-6)}).

\subsection{Example: Linking form of the circle} 
Here we will compute explicitly the linking form of the circle $S^1$. 

As the initial data we have to specify an object $\M$ of a von Neumann category $\CC$
with a nondegenerate
Hermitian form $\psi: \M\to \M^\ast$, $\psi^\dagger = \psi$ and a right unitary action of 
the fundamental group 
$\pi_1(S^1)$ on $\M$. Fix a generator $t\in \pi=\pi_1(S^1)$. Then a unitary
action of $\pi$ on $\M$
is given by a $\CC$-morphism $t: \M\to \M$ such that $t^\ast \psi t=\psi$. 
We will assume that $t-1:\M\to \M$ is injective and has a dense image. 
We will view $\psi$ simply as a scalar product on $\M$. 

Fix a cell decomposition of $K=S^1$, consisting of one zero-dimensional cell $e^0$ and one 
one-dimensional cell $e^1$. The corresponding cell decomposition of the universal covering
$\tilde K$ consists of the lifts $t^ie^0$ and $t^je^1$ of the cells $e^0$ and $e^1$, 
where $i, j \in \Z$, and the boundary homomorphism $d: C_1(\tilde K)\to C_0(\tilde K)$ is given by 
$$d(e^1) = (t-1)e^0.$$
The dual cell decomposition also has only two cells: one zero-dimensional cell $f^0$ and one 
one-dimensional cell $f^1$. The corresponding dual cell decomposition of $\tilde K$ is a shifted
up version of the initial cell decomposition of $\tilde K$ given by $t^ie^0$ and $t^je^1$. 
\begin{figure}[h]
\setlength{\unitlength}{1cm}
\begin{center}
\begin{picture}(10, 4)
\linethickness{0.3mm}
\put(1,1){\vector(1,0){8}}
\put(1,3){\vector(1,0){8}}
\multiput(2.6,1)(2.7,0){3}{\circle*{0.15}}
\multiput(1.3,3)(2.7,0){3}{\circle*{0.15}}
\put(2.6,0.2){$f^0$}
\put(5.25,0.2){$tf^0$}
\put(7.9,0.2){$t^2f^0$}
\put(3.9,0.2){$f^1$}
\put(6.6,0.2){$tf^1$}
\put(1.2,0.2){$t^{-1}f^1$}

\put(1.3,3.2){$e^0$}
\put(3.95,3.2){$te^0$}
\put(6.6,3.2){$t^2e^0$}
\put(2.6,3.2){$e^1$}
\put(5.3,3.2){$te^1$}
\put(8,3.2){$t^{2}e^1$}
\put(0.2,0.9){$\tilde K$}
\put(0.2,2.9){$\tilde K$}
\end{picture}
\end{center}
\end{figure}

As above, the boundary homomorphism of dual chain complex $C'_\ast(\tilde K)$ is given by
$$d(f^1) = (t-1)f^0.$$
Intersection pairing (\ref{(2-1)}) is given by
$$[e^1, f^0] =1, \quad [e^0, f^1] =t^{-1}.$$
Now we have to find a symmetric (i.e. satisfying (\ref{(2-3)})) simplicial approximation $h: C_\ast(\tilde K)\to 
C'_\ast(\tilde K)$ of the identity map $\tilde K\to \tilde K$. Such symmetric simplicial approximation
is given by
$$h(e^0) = \frac{1}{2} (1+t^{-1})f^0, \quad h(e^1) = \frac{1}{2} (1+t^{-1})f^1.$$
One checks that diagram (\ref{(3-18)} in this case has the form
$$
\CD
C_1(\tilde K) @>{t-1}>> C_0(\tilde K)\\
@V{\frac{1}{2}(t+1)}VV  @VV{\frac{1}{2}(t^{-1}+1)}V\\
(C_0(\tilde K))^\ast @>{1-t^{-1}}>> (C_1(\tilde K))^\ast.
\endCD
$$
Therefore, diagram (\ref{(3-19)} becomes
\begin{eqnarray}
\CD
\M @>{t-1}>> \M\\
@V{-\frac{1}{2}(t+1)}VV  @VV{\frac{1}{2}(t^{-1}+1)}V\\
\M @>{t^{-1}-1}>> \M.
\endCD\label{(3-20)}
\end{eqnarray}
This diagram represents the linking form of the circle $S^1$ (cf. \cite{Fa3}, \S 4.1). 
The zero-dimensional extended $L^2$ homology of the circle is $\H_0(S^1;\M)=((t-1):\M\to\M)$,
which is all torsion (since we assume that $t-1:\M\to \M$ has dense image)
and (\ref{(3-20)}) represents a (-1)-Hermitian linking form on $\H_0(S^1;\M)$. 

As the next step we want to express the linking form (\ref{(3-20)}) as a discriminant form of a Hermitian 
form on a projective object, cf. \cite{Fa3}, \S 4. In \cite{Fa3}, \S 4.9 we have described a general 
algorithm for this. Applying it in this situation we arrive at the following. Consider the spectral decomposition
$$2-t -t^{-1} = \int_{0}^{\infty} \lambda dE_\lambda$$
of the self-adjoint non-negative
operator $(t-1)(t-1)^\ast= 2-t -t^{-1}.$ Pick a small $\epsilon>0$ and split
$\M$ as the direct sum
$$\M=\M_\epsilon \oplus Q,\quad\text{where}\quad \M_\epsilon = E_{\epsilon}\M , \quad\text{and}\quad Q=(1 -E_{\epsilon})\M.$$
In other words, $t-1$ is small on $\M_\epsilon$ and it is large on $Q$; therefore, we may make
an excision (cf. \cite{Fa3}, \S 4.6) and cut out $Q$. As the result 
(cf. diagrams (4-14) and (4-15) in \cite{Fa3}) 
we obtain that the linking form
of the circle equals the discriminant form of the following $(-1)$-Hermitian form
\begin{eqnarray}(t-t^{-1})/2: \M_\epsilon \to \M_\epsilon.\end{eqnarray}

In order to make the problem more specific, let $\M=L^2(Z,\mu)$ where $Z$ is a locally compact Hausdorff space with a positive Radon measure $\mu$, cf. \cite{Fa2}, \S 2, Example 7 and also
\cite{Fa3}, \S 7.3. Assume that the unitary 
action of $\pi_1(S^1)$ on $\M=L^2(Z,\mu)$ is given as the 
multiplication on a function $t(z)=\exp(if(z))\in L^\infty_\C(Z,\mu)$, 
where $f:Z\to [-\pi, \pi]$ is a real valued function $f\in L^\infty_\R(Z,\mu)$. 
Then $\M_\epsilon$ is 
$L^2(Z_\epsilon, \mu_\epsilon)$, where 
$$Z_\epsilon =\{z\in Z; |\exp(if(z))-1|^2<\epsilon\},\quad \mu_\epsilon =\mu|_{Z_\epsilon}.$$

The splitting of Theorem 7.7 from \cite{Fa3} produces two torsion objects $\X_+$ and $\X_-$
out of the linking form, $\X_+\oplus \X_-=\H_0(S^1;\M)$. 
Their spectral density functions are given by
$$F_+(\lambda) = \mu\{z\in Z_\epsilon; f(z)\in \sin^{-1}(0,\lambda]\},\quad 
F_-(\lambda) = \mu\{z\in Z_\epsilon; f(z)\in \sin^{-1}([-\lambda,0)\},$$
compare \cite{Fa2}, formula (4-4). Clearly, the above formulae allow to compute (in terms of the function $f$ and the measure $\mu$)
the Novikov-Shubin numbers
of the positive and negative parts $\X_+$, $\X_-$. The result depends only 
on the behavior of function $f$ near points $z\in Z$ where $f(z)$ is "small", 
but in contract with the usual Novikov-Shubin invariants, the points where $f$ is small and positive are treated differently compared to points
where $f$ is small and negative.

The following is a cohomological version of Theorem \ref{thm3.4}.

\begin{theorem}\label{thm3.9} Let $K$ be a closed piecewise 
linear $2q+1$-dimensional 
manifold. Denote by $\pi$ the fundamental group of $K$ and let 
$w:\pi\to \Z_2$ be the first Stiefel - Whitney class of $K$. Consider the 
group ring $\C[\pi]$ with the involution determined by $w$, i.e. 
$\overline g=w(g)g^{-1}$ for $g\in\pi$. Then 
any non-degenerate Hermitian form $\psi:\M\to D_\Lambda(\M)$ in $\Lambda$-mod-$\CC$ and a choice of orientation of $\tilde K$
determine canonically a $(-1)^{q+1}$-Hermitian form 
\begin{eqnarray}
\L: T^{q+1}\to \e(T^{q+1}), \quad \text{where}\quad T^{q+1}=T(\H^{q+1}(K,\M)).
\label{(3-22)}\end{eqnarray}       
\end{theorem} 

We call (\ref{(3-22)}) {\it the cohomological linking form}. 

\subsection{} Given a closed PL manifold $K$, 
the data for the linking form (\ref{(3-6)}) consist of specifying an orientation of $\tilde K$ and
an object $\M$ of a 
Hilbertian von Neumann category $\CC$, a right action of $\pi=\pi_1(K)$ on $\M$ and also a 
Hermitian form $\psi: \M\to \M^\ast$. Geometrically, all these data can be understood as a 
flat bundle $\E\to K$ with fiber $\M$, whose monodromy coincides with the given action of $\pi$
on $\M$, together with a flat bundle map $\langle\, ,\, \rangle : \E\otimes \E\to \C_w$, where
$\C_w$ denotes a flat complex line bundle with monodromy given by the first Stiefel-Whitney class
$w: \pi\to \{1, -1\}$. 

In case when the manifold $K$ is orientable, the bundle $\E$ is simply 
a flat unitary bundle over $K$ with fiber and monodromy in category $\CC$. 

The same data determine the cohomological linking form (\ref{(3-22)}). 

\section{\bf Manifolds with boundary}\label{sec4}

In this section we show how to compute the linking form 
assuming that the given odd-dimensional manifold is a boundary $\partial K$, where $\dim K=2q$.

\subsection{ The intersection form}\label{sec4.1}
Consider a compact oriented $(2q)$-dimensional manifold $K$ with boundary $\partial K$ and a 
flat unitary Hilbertian bundle $\E$ over $K$. The fiber and the monodromy of $\E$ belong to a Hilbertian von Neumann category $\CC$. We have the extended homology
$\H_i(K; \E)$ and $\H_i(\partial K; \E)$. 
We will describe now the intersection form
\begin{eqnarray}\I: P_q(K) \to (P_q(K))^\ast,\label{(4-1)}\end{eqnarray}
where we denote
\begin{eqnarray}P_q(K) = P(\H_q(K;\E)/\ker[j:P(\H_q(K;\E))\to P(\H_q(K,\partial K;\E))].\label{(4-2)}\end{eqnarray}
Recall that $P(\H_q(K;\E))$ 
denotes the projective part of $\H_q(K;\E)$, the extended $L^2$-homolo\-gy with coefficients in $\E$.
The intersection form (\ref{(4-1)}) is defined as follows. We have the following commutative diagram
$$
\CD
\H_q(K;\E) @>j>> \H_q(K, \partial K;\E)\\
@V{\simeq}VV @VV{\simeq}V\\
\H^q(K,\partial K;\E) @>j^\ast>> \H^q(K;\E)
\endCD
$$
in the extended abelian category $\E(\CC)$, where the vertical maps are Poincar\'e duality isomorphisms and the horizontal maps are induced by the inclusion $j$. 
Applying the functor of projective part we obtain the commutative diagram
$$
\CD
P(\H_q(K;\E)) @>j>> P(\H_q(K, \partial K;\E))\\
@V{\simeq}VV @VV{\simeq}V\\
P(\H^q(K,\partial K;\E)) @>j^\ast>> P(\H^q(K;\E)).
\endCD
$$
By the Universal Coefficients
Theorem, $P(\H^q(K,\partial K;\E))\simeq P(\H_q(K, \partial K;\E))^\ast$, and also
$P(\H^q(K;\E))\simeq P(\H_q(K;\E))^\ast$; moreover, the homomorphism 
$j^\ast$ (downstairs) is dual to $j$ (upstairs). The morphism $P(\H_q(K;\E))\to P(\H^q(K;\E))\simeq
(P(\H_q(K;\E))^\ast$ (acting from the left upper corner to the right lower corner)
vanishes on $\ker j$ and takes the values in $((P(\H_q(K;\E))/\ker j)^\ast\subset 
(P(\H_q(K;\E))^\ast$.
Recall that by the definition, $P(\H_q(K;\E))/\ker j = P_q(K)$. Hence the above diagram
determines  
$$\I: P_q(K)\to (P_q(K))^\ast,$$ 
which by definition is the intersection form. 

The intersection
form $\I$ is $(-1)^q$-Hermitian and is {\it weakly nondegenerate},
i.e. the morphism (\ref{(4-1)}) has zero kernel. Note that in general the intersection form $\I$ {\it
fails to be non-degenerate}.
In fact, our goal in this section is to compute the discriminant
form of $\I$
(cf. \cite{Fa3}, \S 4.4), which describes the way $\I$ degenerates.

The intersection form $\I$ is nondegenerate if $K$ is closed $\partial K=\emptyset$.

As we know from Theorem 3.4, in the situation above  
there also exists the linking form
\begin{eqnarray}
\L : T_{q-1}(\partial K) \to \e(T_{q-1}(\partial K)),\label{(4-3)}
\end{eqnarray}
where
$T_{q-1}(\partial K)$ denotes $T(\H_{q-1}(\partial K; \E))$, 
the torsion part of the extended homology. The linking form 
$\L$ is non-degenerate and $(-1)^{q}$-Hermitian.

The following is the main result of this section.

\begin{theorem}\label{thm4.2}
Let $\X \subset T_{q-1}(\partial K)$ denote the image
of the boundary homomorphism
\begin{eqnarray}
\X = \im[\partial: T(\H_{q}(K,\partial K;\E)) \to T_{q-1}(\partial K)].
\label{(4-4)}\end{eqnarray}
Then $\X$ is contained in its annihilator $\X^{\perp}$ with respect to the linking form
(\ref{(4-3)}) and the induced by $\L$ the non-degenerate 
form on $\X^{\perp}/\X$ is congruent to the discriminant form
(cf. \cite{Fa3}, \S 4.4) of the intersection form $\I$, cf (\ref{(4-1)}).
\end{theorem}

Let's explain the terms used in the statement of Theorem \ref{thm4.2}. Suppose that
 $\L: T\to \e(T)$ is a torsion Hermitian form. If $i: \X\to T$ is an embedding (in the sense of the 
abelian category $\eca$) then the annihilator $\X^\perp\subset T$ is defined as the kernel of 
$\e(i)\circ \L: T\to \e(\X)$. The assumption that $\X$ is contained in $\X^\perp$ is equivalent to 
vanishing of the composite $\X@>i>> T@>{\e(i)\circ \L}>>\e(\X)$. The induced torsion Hermitian
form on $\X^\perp/\X$ is defined as follows. Consider the commutative diagram
$$
\CD
\X@>i>>T@>{\e(i)\circ \L}>>\e(\X)\\
@V{1}VV @V{\L}VV @V{1}VV\\
\X@>>{\L\circ i}> \e(\T)@>>{\e(i)}>\e(\X).
\endCD
$$
Both rows are chain complexes and the homology of the upper row is $\X^\perp/\X$. Homology
of the lower row is canonically isomorphic to $\e(\X^\perp/\X)$ (use Lemma in \cite{Fa3}, \S 1.1).
Since the vertical arrows in the above diagram are isomorphisms, it 
determines an isomorphism $\X^\perp/\X\to \e(\X^\perp/\X)$, i.e. a nondegenerate form on $\X^\perp/\X$.

\begin{corollary} The intersection form $\I$ is non-degenerate
if and only if the torsion submodule $\X \subset T_{q-1}(\partial K)$, given by (\ref{(4-4)}), is
a metabolizer of the linking form $\L$
(i.e. $\X = \X^\perp$). \qed
\end{corollary}

\begin{corollary} Suppose that under the conditions of Theorem \ref{thm4.2} 
it is additionally known that $T(\H_{q-1}(K;\E)) =0$. Then the
linking form $\L: T_{q-1}(\partial K) \to \e(T_{q-1}(\partial K))$ 
is congruent to the discriminant of the intersection form $\I$. 
\end{corollary} 
\begin{proof} If $T(\H_{q-1}(K;\E)) =0$ then by the Poincar\'e duality 
$T(\H_{q}(K,\partial K;\E)) =0$, and hence $\X=0$ (cf. (4.4)). Therefore in this case $\X^\perp$
coincides with $T_{q-1}(\partial K)$. The result now follows from Theorem \ref{thm4.2}.
 \end{proof}

In the proof of Theorem \ref{thm4.2} we will use the following Lemma.

\begin{lemma}\label{l4.5} Suppose that
$\H_\ast = (\dots \to \H_{i+1} \to \H_i \to \H_{i-1} \to \dots) $
is an exact sequence in the extended abelian category $\eca$. Form two chain complexes
$$P_\ast = (\dots \to P(\H_{i+1}) \to P(\H_i) \to P(\H_{i-1}) \to \dots)$$
and
$$T_\ast = (\dots \to T(\H_{i+1}) \to T(\H_i) \to T(\H_{i-1}) \to \dots)$$
by taking the projective and torsion parts of $\H_\ast$
respectively. Consider homology of $P_\ast$ and $T_\ast$ in $\eca$, which we denote
$\H_i(P_\ast)$ and $\H_{i}(T_\ast)$ correspondingly.
Then for any $i$ there is a natural isomorphism
$$\H_i(P_\ast) \to \H_{i-1}(T_\ast).$$
\end{lemma}
\begin{proof} We have a short exact sequence $0\to T_\ast \to \H_\ast \to P_\ast \to 0$
of chain complexes in $\eca$.
The corresponding long homological sequence provides the required isomorphism, since
$\H_\ast$ is exact. 
\end{proof}

\subsection{Proof of Theorem \ref{thm4.2}} Consider the following commutative diagram
with exact rows
\begin{eqnarray}
\CD
\dots\to \H_q(K,\partial K;\E) @>{\partial}>> \H_{q-1}(\partial K;\E) @>{i_\ast}>> \H_{q-1}(K;\E)\to \dots\\
@V{\simeq}VV   @V{\simeq}VV  @V{\simeq}VV\\
\dots \to  \H^q(K;\E) @>{i^\ast}>>  \H^{q}(\partial K;\E) @>{\delta}>> 
\H^{q+1}(K,\partial K; \E)\to \dots
\endCD\label{(4-5)}
\end{eqnarray}
The vertical arrows denote the isomorphisms of Poincare duality. 
We will apply to this diagram the functors of torsion and projective 
parts and then use the natural isomorphisms of Lemma \ref{l4.5}. Applying the functor of the
torsion part we obtain the commutative diagram
\begin{eqnarray}
\CD
T_q(K,\partial K) @>{\partial}>> T_{q-1}(\partial K) @>{i_\ast}>> T_{q-1}(K)\\
@V{\simeq}VV   @V{\simeq}VV  @V{\simeq}VV\\
\e(T_{q-1}(K)) @>{\e(i_\ast)}>>  \e(T_{q-1}(\partial K)) @>{\e(\partial)}>>
\e(T_{q}(K,\partial K)).
\endCD\label{(4-6)}
\end{eqnarray}
The vertical morphisms here can be viewed as three different linking forms 
induced by the Poincare duality. 
The homology of the upper horizontal sequence in (\ref{(4-6)})
can be easily identified with $\X^{\perp}/\X$. The lower horizontal sequence is dual
to the upper one.The middle vertical isomorphism of (\ref{(4-6)}) induces an isomorphism from
the homology of the upper horizontal sequence to the homology of the lower one;
this isomorphism is by the definition 
the induced form $\X^{\perp}/\X \to \e(\X^{\perp}/\X)$.

Now we apply the functor of projective part $P$ to diagram (\ref{(4-5)}) (at the place shifted by one)
to get the commutative diagram
\begin{eqnarray}
\CD
P(\H_q(K;\E)) @>j>> P(\H_q(K,\partial K;\E)) @>{\partial}>> P(\H_{q-1}(\partial K;\E))\\
@V{\simeq}VV   @V{\simeq}VV  @V{\simeq}VV\\
P(\H_q(K,\partial K;\E))^\ast @>{j^\ast}>>  P(\H_q(K;\E))^\ast @>>>
P(\H_{q}(\partial K;\E))^\ast .
\endCD\label{(4-7)}
\end{eqnarray}
The left square of this diagram appeared above in subsection \ref{sec4.1} 
in the definition of the intersection form $\I$.
We want to compute the homology of the two horizontal sequences in the middle and to identify the
results with the corresponding homology of the rows of (\ref{(4-6)}) using Lemma \ref{l4.5}. We find that
the diagram
\begin{eqnarray}
\CD
P_q(M) @>\I>> P_q(M)^\ast @>>> 0\\
@V1VV @V1VV @VVV\\
P_q(M) @>\I>> P_q(M)^\ast @>>> 0,
\endCD
\label{(4-8)}
\end{eqnarray}
which is built entirely out of the intersection pairing, has homology of the horizontal
sequences in the middle term identical (in a canonical way) to the homology of (\ref{(4-7)}).
Thus we obtain that
(\ref{(4-8)}) can be mapped onto (\ref{(4-7)}) 
inducing isomorphisms of the extended homology in the
middle places. Comparing with the definition of discriminant form in \cite{Fa3}, \S 4.4, we see
that this isomorphism between $\X^\perp /\X$ and $(\I: P_q(M)\to P_q(M)^\ast)$ gives a
congruence between the induced form on $\X^\perp /\X$ and the discriminant form of the
intersection form. \qed

\section{\bf  Hyperbolicity of the linking form}\label{sec5}

\begin{theorem}\label{thm5.1}
Suppose that $M$ is a closed orientable
$(2q)$-dimensional manifold and
$\E \to M$ is a flat unitary Hilbertian bundle having the fiber and the monodromy in a finite von
Neumann category
$\CC$. Suppose that $W \subset M$ is a closed
codimension one submanifold, which separates $M$ into two parts $M_+$ and $M_-$. If the
middle dimensional extended homology $\H_q(M;\E)=0$ vanishes, then the linking form of $W$
$$\L : T(\H_{q-1}(W;\E)) \to \e(T(\H_{q-1}(W;\E)))$$
is hyperbolic (cf. \cite{Fa3}, \S 5). \end{theorem}

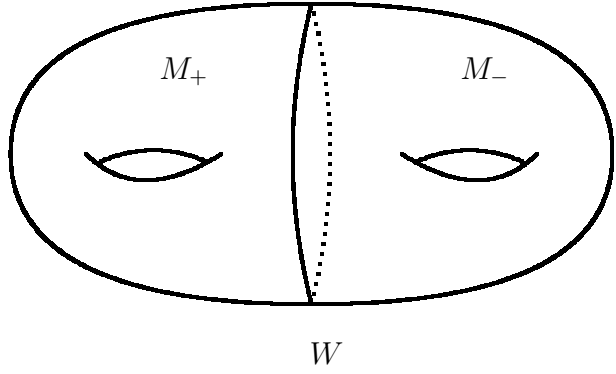
\begin{figure}[h]
\setlength{\unitlength}{1cm}
\begin{center}
\begin{picture}(10, 6)
\linethickness{0.4mm}
\qbezier(5,1)(1,1)(1,3)
\qbezier(5,5)(1,5)(1,3)
\qbezier(5,5)(9,5)(9,3)
\qbezier(5,1)(9,1)(9,3)
\qbezier(5,1)(4.5, 3)(5,5)
\qbezier[30](5,1)(5.5, 3)(5,5)
\qbezier(2,3)(2.7,2.3)(3.8,3)
\qbezier(8,3)(7.3,2.3)(6.2,3)
\qbezier(2.2,2.9)(2.9,3.2)(3.6,2.9)
\qbezier(7.8,2.9)(7.1,3.2)(6.4,2.9)
\put(5,0.2){$W$}
\put(3,4){$M_+$}
\put(7,4){$M_-$}
\end{picture}
\end{center}
\end{figure}

\begin{proof} Applying Theorem \ref{thm4.2} to two manifolds $M_+$ and $M_-$,
we obtain two submodules $\X_\pm \subset T(\H_{q-1}(W;\E))$, described in Theorem 
\ref{thm4.2}.

First, we observe that each of $\X_\pm$ is a metabolizer, i.e. $\X_\pm ^\perp = \X_\pm$. Indeed,
by Theorem \ref{thm4.2}, the factor $\X_\pm ^\perp / \X_\pm$ is isomorphic to the cokernel of the
intersection pairing $\I_\pm : P_q(M_\pm) \to (P_q(M_\pm))^\ast$ (we use the notations introduced
in section 4.1). However, our assumption $\H_q(M;\E) = 0$ implies that $P_q(M_\pm) =0$.

Now we show that the metabolizers $\X_+$ and $\X_-$ are "disjoint", i.e.
$\X_+ \cap \X_- =0$ (here the intersection is understood as intersection of two subobjects
of an object of an abelian category).
Indeed, according to Theorem \ref{thm4.2} and the previous arguments, we may
identify
$\X_\pm$ with the kernel of the morphism $T(\H_{q-1}(W;\E)) \to T(\H_{q-1}(M_\pm;\E))$
induced by the inclusion. Using the Mayer - Vietoris sequence
\begin{eqnarray}
0 \to \H_{q-1}(W;\E) \stackrel{{i_+}_\ast + {i_-}_\ast}\longrightarrow
\H_{q-1}(M_+;\E)\oplus \H_{q-1}(M_-;\E)\to \H_{q-1}(M;\E)
\label{(5-1)}
\end{eqnarray}
and our assumption
$\H_q(M;\E) = 0$, we obtain $\X_+ \cap \X_- =0$.

Now we want to show that $\X_+ + \X_- = T(\H_{q-1}(W;\E))$. 
We examine the Mayer - Vietoris
sequence (\ref{(5-1)}) again. Using the Poincar\'e duality we find that our assumption
$\H_q(M;\E) = 0$ implies $T(\H_{q-1}(M;\E)) = 0$. Therefore we conclude that the direct sum
$T(\H_{q-1}(M_+;\E))\oplus T(\H_{q-1}(M_-;\E))$ belongs to the image of ${i_+}_\ast + {i_-}_\ast$.
Thus we get an isomorphism
$${i_+}_\ast + {i_-}_\ast : T(\H_{q-1}(W;\E))\stackrel{\simeq}{\longrightarrow} 
T(\H_{q-1}(M_+;\E))\oplus T(\H_{q-1}(M_-;\E)),$$
which implies $\X_+ + \X_- = T(\H_{q-1}(W;\E))$.

As a result, $\X_\pm$ are mutually complementary metabolizers  in  
$T(\H_{q-1}(W;\E)))$, and this implies hyperbolicity of the linking form of $W$. 
\end{proof}

\section{\bf Novikov - Shubin signatures}\label{sec6}

In this section we define a new numerical invariant of torsion Hermitian forms which we call the Novikov - Shubin signature.

\subsection{} Let $K$ be a closed oriented $2q+1$-dimensional manifold. Let $\E\to K$ be a
flat unitary bundle with fiber and monodromy in a superfinite von Neumann category $\CC$, cf. \cite{Fa3}, \S 7. Then we have the $(-1)^{q+1}$-Hermitian linking form 
\begin{eqnarray}\L: T_q\to \e(T_q).\label{(6-1)}
\end{eqnarray}
Here $T_q$ denotes the torsion part of the extended $L^2$-homology, $T_q=T(\H_q(K;\E))$.
The form 
\begin{eqnarray}
\cases 
\L, \quad\text{if}\quad q \quad\text{is odd},\\
-\sqrt{-1}\, \L, \quad\text{if}\quad q \quad\text{is even}
\endcases
\label{(6-2)}
\end{eqnarray}
is Hermitian, and by Theorem 7.7 from \cite{Fa3}, $T_q$ splits canonically as a direct sum
\begin{eqnarray*}T_q = (T_q)_+ \oplus (T_q)_-\end{eqnarray*}
of a positive and a negative definite parts of (\ref{(6-2)}). 
This splitting contains all the information about 
the linking form $\L$, cf. \cite{Fa3}, \S 7.

Hence numerical invariants of the linking form $\L$ can be obtained by computing numerical invariants
of the positive $(T_q)_+$ and negative $(T_q)_-$ parts of $T_q$. 

\subsection{} The best known invariant of torsion objects is 
{\it the Novikov-Shubin number},
introduced by S.P. Novikov and M.A. Shubin in \cite{NS}. In the context of von Neumann categories
the Novikov - Shubin number was studied in \cite{Fa1}, \cite{Fa2}.
In order to define the Novikov - Shubin number one should have specified 
(additionally to the above mentioned data) a trace
on the category $\CC$, cf. \cite{Fa2}, \S 3.9.  
The Novikov-Shubin number of a torsion object $\X$
is denoted $\ns(\X)$, or $\ns_{\tr}(\X)$, if we want to emphasize the used trace $\tr$.

It was observed in \cite{Fa1}, that a more convenient invariant (which is equivalent to $\ns(\X)$)
is given by
\begin{eqnarray*}\c(\X) = \ns(\X)^{-1}\in [0, \infty].\end{eqnarray*}
It is 
called {\it the Novikov - Shubin capacity of} $\X$. 
The advantage of capacity against the Novikov-Shubin number
consists in the fact that it adequately describes the size of a torsion 
object: larger torsion objects have larger capacity, trivial object 
has zero capacity, there may also exist torsion objects having infinite
capacity.

\subsection{}
We my use the Novikov-Shubin capacity to obtain invariants of linking forms by setting
\begin{eqnarray}
\c_+(\L) = \c((T_q)_+), \quad \c_-(\L) = \c((T_q)_-).
\label{(6-5)}
\end{eqnarray}
The pair of numbers $(\c_+(\L), \c_-(\L))$ will be called {\it the Novikov - Shubin signature of } $\L$.

In fact, 
\begin{eqnarray}
\max\{\c_+(\L), \c_-(\L)\} = \c(T_q) = \c(T(\H_q(K;\E)))\label{(6-6)}
\end{eqnarray}
(cf. \cite{Fa2}, (3-22)). Therefore, only one of the numbers $\c_+(\L)$ or $\c_-(\L)$
may actually carry a new information (in case it is smaller than the Novikov - Shubin capacity
$\c(T_q)$).

Invariants (\ref{(6-5)}) may be useful in order to show that the linking form $\L$ 
is not hyperbolic.

\begin{proposition}\label{prop6.4}
(a) If linking form (\ref{(6-1)}) is hyperbolic then 
$$\c_+(\L) = \c_-(\L).$$
(b) If $\L = \L_1\perp \L_2$ is an orthogonal sum of forms $\L_1$ and $\L_2$, then 
$$\c_\pm(\L) =\max\{\c_\pm(\L_1), \c_\pm(\L_2)\}.$$
\end{proposition} 
\begin{proof} (a) follows from \cite{Fa3}, Theorem 7.10. (b) follows from Proposition 3.10 in \cite{Fa2}.\end{proof}

\subsection{Example}\label{sec6.5} 
We will compute now the invariants (\ref{(6-5)}) in the case of manifolds
with $\pi_1(K)=\Z$. We will assume here that the unitary flat bundle $\E$ has the fiber $\ell^2(\Z)$
and the monodromy is given by the standard
action of  $\Z$ on $\ell^2(\Z)$. 
It is well known (using Fourier transform) that we may equivalently think of $\ell^2(\Z)$ as of 
$\ell^2(S^1)$, such that the action of the generator $t\in \Z$ is given as the multiplication 
$\ell^2(S^1)\to \ell^2(S^1)$ by the identity function $t=\id : S^1\to S^1\subset \C$. 
The von Neumann category $\CC$
generated by the $\C[\Z]$-module $\ell^2(\Z)$ is isomorphic to the category of finite dimensional
measurable fields of Hilbert spaces over $S^1$, which is superfinite, according to \cite{Fa3}, \S 7.3.
A trace on this category can be specified by choosing a measure on the circle $S^1$. (Note that 
there are traces on this category which do not come from measures on $S^1$).

We will fix this measure as follows. Let $J\subset S^1$ be a closed 
interval and let $\mu_J$ be the measure on $S^1$, which is the usual Lesbegue measure restricted to $J$ and which is 
zero outside $J$. We will assume that $J$ is small enough, cf. below. 

Let $K$ be a $(2q+1)$-dimensional closed oriented manifold with $\pi_1(\Z)$.
Consider the universal covering $\tilde K\to K$. It has $\Z$ as the group of covering transformations. The homology $H_q(\tilde K;\C)$ is a direct sum of a free $\C[\Z]$-module and of a torsion 
$\C[\Z]$-module. The torsion submodule also splits as a finite 
direct sum of modules of the form
$$M_{c,m}=\C[\Z]/(t-c)^m\C[Z]$$ 
where $c\in \C$ and $m\ge 0$, $m\in \Z$. 
We will denote by $\mathcal S_q\subset \C$ the set of 
those $c$, which lie on the unit circle $S^1$. $\mathcal S_q(K)$ is a finite set, 
which we will call {\it the support of}
$K$. 

We will assume that interval $J\subset S^1$ is so small that
it contains at most one point $c\in \mathcal S_q(K)$.

For $c\in \mathcal S_q(K)\cap J$ consider the $(t-c)$-adic part $(H_q(\tilde K;\C))_c$ of $H_q(\tilde K;\C)$, 
i.e. the set of all homology classes,
which can be killed by a power $(t-c)^m$ with some large $m$. 
There is a $(-1)^{q+1}$-Hermitian {\it Blanchfield pairing}
\begin{eqnarray}
\{\, ,\, \}: (H_q(\tilde K;\C))_c\times (H_q(\tilde K;\C))_c\to \RR/\C[\Z]_c,
\label{(6-7)}
\end{eqnarray}
where $\RR$ denotes the field of fractions of $\C[\Z]$ and $\C[\Z]_c$ denotes the ring
rational function in $t$ such that their denominators are prime to $t-c$, cf. \cite{FL}, \S 9.
The involution in $\RR/\C[\Z]_c$ is given by the complex conjugation and by
$\overline t=t^{-1}$.

We will explain below (cf. \ref{sec6.8})
how the Blanchfield pairing (\ref{(6-7)}) determines nonnegative integers
$n_j^+(c)$ and $n_j^-(c)$, where $j=1, 2, \dots$ such that only finitely many of them are nonzero.

Now we define {\it odd and even height numbers} as follows:
\begin{eqnarray}
h_{\text{odd}}(c) = \max\{j\, \, \text{odd, such that}\,\,  n_j^+(c)>0\, \, \text{or}\, \, n_j^-(c)>0\},\label{(6-8)}\end{eqnarray}
\begin{eqnarray}h_{\text{ev}}^\pm(c) = \max\{j\, \, \text{even, such that}\, \, n_j^\pm(c) >0\}.\label{(6-9)}
\end{eqnarray} 

\begin{theorem}\label{thm6.6} In the situation described in 6.5, assume that 
$J\cap \mathcal S_q(K)$ consists of a single interior point $c\in J$. 
Then the positive and negative 
Novikov-Shubin capacities (\ref{(6-5)}) of the linking form
\begin{eqnarray}
\L: T_q\to \e(T_q), \quad T_q=T(\H_q(K;\ell^2(\Z)))\label{link}
\end{eqnarray}
with respect to the trace determined by the measure $\mu_J$ are given by
\begin{eqnarray}\c_+(\L) = \max\{h_{\text{odd}}(c), h^+_{\text{ev}}(c)\}, \quad 
\c_-(\L) = \max\{h_{\text{odd}}(c), h^-_{\text{ev}}(c)\}.\end{eqnarray}
\end{theorem}
Note that 
$$\max\{h_{\text{odd}}(c), h^+_{\text{ev}}(c), h^-_{\text{ev}}(c)\}$$
coincides with the height
of $(H_q(\tilde K;\C))_c$, and hence by  \cite{Fa2}, Theorem 4.13, it coincides with the Novikov-Shubin capacity of $\H_q(K;\ell^2(\Z))$ with respect to $\mu_J$.

Theorem \ref{thm6.6} demonstrates that the numbers $\c_+(\L)$ and $\c_-(\L)$ 
can be arbitrary except that their maximum
is fixed by relation (\ref{(6-6)}). 

Let us consider now situation when $\mathcal S_q(K)\cap J$ is a single point, 
which is an end point of the interval $J$. We assume
that the circle $S^1$ is canonically oriented (anticlockwise) and hence we may speak about the
initial and the terminal points of $J$. 

\begin{theorem}\label{thm6.7} In the situation described in \ref{sec6.5}, assume that the intersection
$J\cap \mathcal S_q(K)$ consists of a single point $c$, which is the terminal point of the interval $J$. 
Then the positive and negative 
Novikov-Shubin capacities (\ref{(6-5)}) of the linking form
$$\L: T_q\to \e(T_q), \quad T_q=T(\H_q(K;\ell^2(\Z)))$$
with respect to the trace determined by the measure $\mu_J$ are given by
\begin{eqnarray}
\c_+(\L) = \max\{h^+_{\text{ev}}(c), h^-_{\text{odd}}(c)\}, \quad 
\c_-(\L) = \max\{h^-_{\text{ev}}(c), h^+_{\text{odd}}(c)\}.
\label{(6-12)}
\end{eqnarray}
\end{theorem}

\subsection{Numbers $n_j^\pm(c)$}\label{sec6.8} 
Here we will explain how the Blanchfield form (\ref{(6-7)}) determines
the numbers $n_j^\pm$, where $j=1, 2, \dots$, which appear in Theorems (\ref{(6-8)}), (\ref{(6-9)}) and (\ref{(6-12)}).
Here we essentially follow \cite{FL}, \S 2. 

Any element $f\in \RR/\C[\Z]_c$ can be uniquely represented in the form
$$f=\alpha_1g+\alpha_2g^2+\alpha_3g^3 +\dots,\quad\text{where}\quad 
g=ic(t-c)^{-1},$$
with $\alpha_j\in \C$;
only finitely many $\alpha_j$'s are nonzero. 
Note that  $g\in \RR/\C[\Z]_c$ is "real", i.e. $\overline g =g$ assuming that $c\in S^1$.

Given 
$x, y\in (H_q(\tilde K;\C))_c$, the value $\{x,y\}\in \RR/\C[\Z]_c$ of the Blanchfield form (\ref{(6-7)}) 
can be uniquely expressed as 
\begin{eqnarray*}
\{x,y\} = \alpha_1(x,y)g+\alpha_2(x,y)g^{2}+\alpha_3(x,y)g^{3}+\dots ,
\end{eqnarray*}
and this defines a sequence of $(-1)^{q+1}$-Hermitian forms
\begin{eqnarray*}
\alpha_j: (H_q(\tilde K;\C))_c\times (H_q(\tilde K;\C))_c\to \C, \quad j=1, 2, \dots.
\end{eqnarray*}
Let $T_j\subset (H_q(\tilde K;\C))_c$ denote the subspace, 
consisting of cycles $z\in (H_q(\tilde K;\C))_c$,
such that $(t-c)^jz=0$. 
Now, we define the numbers $n_j^+(c)$ and $n_j^-(c)$ as the numbers of positive and negative
squares (correspondingly) in the diagonal representation of the Hermitian form 
\begin{eqnarray*}
\cases
\alpha_j: T_j\times T_j\to \C,\quad\text{if $q$ is odd},\\
-\sqrt{-1}\alpha_j: T_j\times T_j\to \C,\quad\text{if $q$ is even}.
\endcases
\end{eqnarray*}
\subsection{Beginning of the proof of Theorems \ref{thm6.6} and \ref{thm6.7}} 
As in \cite{FL}, let $\OO$ denote the ring of germs at the origin of complex valued holomorphic functions $f: (-\epsilon, \epsilon)\to \C$. An element
of $\OO$ can also be represented by a power series
$$f(\tau) =\sum_{n\ge 0} a_n\tau^n, \quad a_n\in \C$$
having a nonzero radius of convergence. The ring operations are given by pointwise addition and multiplication. The involution in $\OO$ is given by $\overline \tau =\tau$ and by the 
complex conjugation. 

The algebraic structure of $\OO$ is extremely simple since 
it is a discrete valuation ring. Its maximal
ideal $\m\subset \OO$ is given by $\m=\{f\in \OO; f(0)=0\}$. Any finitely generated $\OO$-module
is a direct sum of a free module and a torsion submodule. The torsion submodule can be represented as a direct sum of finitely many modules of the form $\OO/\m^k$, where $k\in \Z$.

Let $\CC'$ denote the category of finitely generated $\OO$-modules and let 
$\T'\subset \CC'$ be the full subcategory generated by torsion modules. 

We will describe a duality
$(\e', s')$ in $\T'$, cf. \cite{Fa3}, \S 1. The functor $\e':\T'\to \T'$ is given by 
$\e'(X)=\overline{\Hom}_{\OO}(X,\M/\OO)$, where $X\in \ob(\T')$, and 
$\M$ denotes the ring of fractions of $\OO$. The bar over $\Hom$ means that we consider the set
of all anti-homomorphisms, compare \cite{Fa3}, \S 1.4.  
Note that 
$\M$ can be identified with the ring of germs of meromorphic function $f:(-\epsilon, \epsilon)\to \C$.
The canonical isomorphism $s': X\to \e'\e'(X)$ is given by formula (1-13) in \cite{Fa3} (evaluation and
conjugation).

As in \ref{sec6.5}, $\CC$ denotes the von Neumann 
category of measurable fields of finite dimensional Hilbert spaces over $S^1$ and $\T\subset \eca$ 
denotes its torsion subcategory. Fix a point $c\in S^1$.
We will describe now a covariant functor 
\begin{eqnarray*}
F_c: \T' \to \T,\end{eqnarray*}
compatible with the dualities $\e$ in $\T$ (cf. \cite{Fa3}, \S 3) and $\e'$ in $\T'$. 
Given a torsion object
$X\in \ob(\T')$, we may find a free resolution
\begin{eqnarray}
0\to \OO^n \stackrel{\alpha}\to \OO^n \to X\to 0,
\label{(6-17)}\end{eqnarray}
where $\alpha=(a_{ij}(\tau))$ is an $n\times n$-matrix over $\OO$, such that $\alpha$ 
is invertible over the ring of fractions $\M$. It follows that 
for some small $\epsilon>0$ all functions $a_{ij}(\tau)$ are defined
and analytic for $\tau\in [-\epsilon, \epsilon]$ and the matrix $(a_{ij}(\tau))$ is invertible for
$\tau\in [-\epsilon, \epsilon]$,  $\tau\ne 0$. Let $I_\epsilon\subset S^1$ be the interval
$I_\epsilon=\{c\exp(i\phi); \phi\in (-\epsilon, \epsilon)\}$ on the unit circle. We define
\begin{eqnarray}
F_c(X) = ((a_{ij}): \ell^2(I_\epsilon,\mu)^n\to \ell^2(I_\epsilon,\mu)^n)\, \in \ob(\T).
\label{(6-18)}\end{eqnarray}
Here $\mu$ denotes the Lesbegue measure on $S^1$ 
and $(a_{ij})$ acts by the matrix multiplication. $\ell^2(I_\epsilon,\mu)^n$ is the space of $L^2$ section
of a trivial bundle of rank $n$ over $I_\epsilon$.
It is clear that the obtained object $F_c(X)$
is torsion and it is
independent of the choice of $\epsilon$. To show that $F_c$ is functorial, suppose that 
\begin{eqnarray*}
0\to \OO^{n'} \stackrel{\alpha'}\to \OO^{n'} \to X'\to 0\end{eqnarray*}
is another torsion object of $\T'$. Any $\OO$-homomorphism $h: X\to X'$ leads to a commutative
diagram
\begin{eqnarray*}
\CD
0 @>>> \OO^n @>{\alpha}>> \OO^n @>>> X @>>> 0,\\
 &  & @V{\phi}VV @V{\psi}VV @VVhV  \\
0 @>>> \OO^{n'} @>{\alpha'}>> \OO^{n'} @>>> X' @>>> 0,
\endCD\end{eqnarray*}
and hence we obtain the following morphism of $\T'$
\begin{eqnarray*}
\CD
(\alpha: & (\ell^2(I_\delta, \mu))^n @>>> (\ell^2(I_\delta, \mu))^n)\\
& @V{\phi}VV @V{\psi}VV \\
(\alpha': & (\ell^2(I_\delta, \mu))^{n'} @>>> (\ell^2(I_\delta, \mu))^{n'}),
\endCD\end{eqnarray*}
which, viewed as morphism of the extended category $\eca$ (cf. \cite{Fa2}, \S 1.3),
clearly depends only on $h$ and does not depend on the choice of $\phi$ and $\psi$; 
we will denote it $F_c(h): F_c(X)\to F_c(X')$. Here $\delta>0$
is some small number. 

Let us show that there is a canonical isomorphism $\gamma: F_c\circ \e' \to \e\circ F_c$.
If $\X\in \ob(\T')$ is given by 
resolution (\ref{(6-17)}), where $\alpha =(a_{ij})$, then it is clear that $\e'(X)$ has resolution
\begin{eqnarray*}
0\to \OO^n \stackrel{\alpha^\ast}\longrightarrow \OO^n \to \e'(X)\to 0,
\end{eqnarray*}
where $\alpha^\ast=(\overline{a_{ji}})$ (transposition and conjugation) and hence $F_c(\e'(X))$
is given by
\begin{eqnarray}
((\overline{a_{ji}}): \ell^2(I_\epsilon,\mu)^n\to \ell^2(I_\epsilon,\mu)^n)\, \in \ob(\T).
\label{(6-23)}
\end{eqnarray}
The same object (\ref{(6-23)}) 
is obtained if we first apply $F_c(X)$ (cf. (\ref{(6-18)})) and then the duality $\e(F_c(X))$.

Using $F_c$ and $\gamma$ we may apply the construction of push forward Hermitian forms, 
cf. \S \ref{sec1.1}. Namely, if $\psi: X\to \e'(X)$ is a Hermitian form in $\T'$ then the push forward
form $\phi: F_c(X)\to \e(F_c(X))$ is a Hermitian form in $\T$. 

\begin{theorem}\label{claim}
Let $c\in S^1$ and let $\pi_1(K)\to \OO$ be given by 
$t\mapsto (\tau\mapsto c\exp(i\tau))\in \OO$. This determines a local system $\E^\OO$ of 
free $\OO$-modules over $K$ of rank one. The homology $H_q(K;\E^\OO)$ is a direct sum of a free
module and a torsion submodule; we will denote the latter by $T_q^\OO$. As in \cite{FL}, we have
the linking form $\L^\OO: T_q^\OO\to \e'(T_q^\OO)$. On the other hand, consider an interval
$I_\epsilon\subset S^1$, containing $c$ in its interior, and let $\CC_\epsilon$ be the von Neumann category of measurable fields of finite dimensional Hilbert spaces over $I_\epsilon$. 
The multiplication operator $\ell^2(I_\epsilon, \mu)\to \ell^2(I_\epsilon, \mu)$
by the function $\id: I_\epsilon \to S^1\subset \C$ determines a flat unitary bundle $\E$ over $K$ 
with fiber and monodromy in category $\CC_\epsilon$. Let 
$T_q = T(\H_q(K;\E) = T(\H_q(K; \ell^2(I_\epsilon, \mu)))$ denote the
torsion part of the extended $L^2$ homology in dimension $q$.
Then for small
enough $\epsilon$ the linking form $\L^\epsilon: T_q\to \e(T_q)$ is congruent to the push forward, 
with respect to the pair $(F_c, \gamma)$, cf. \S 1.1,
of the form $\L^\OO$ .\end{theorem}

Note that $\CC_\epsilon$ can be naturally viewed as a full subcategory of $\CC$; 
hence the linking form $\L^\epsilon$
can be viewed as a Hermitian form in $\T$. 

$\L^\epsilon$ is a localized version (at point $c$) of the form (\ref{link}).

Theorem \ref{claim} allows to recover the form $\L^\epsilon$ 
in terms of the linking form of the deformation, 
which was studied in \cite{FL}. 

\subsection{Proof of Theorem \ref{claim}} Consider two mutually dual triangulations of $K$. 
Consider diagram (\ref{(3-18)}.
Tensoring it by $\OO$ we obtain a commutative diagram of the form
\begin{eqnarray}
\CD
F_1@>{\alpha}>> F_0\\
@V{(-1)^{q+1}h^\ast}VV @VVhV\\
F^\ast_0@>{\alpha^\ast}>> F_1^\ast
\endCD\label{(6-24)}
\end{eqnarray}
consisting of free finitely generated $\OO$-modules such that: 
\begin{itemize}
\item[(a)] the torsion submodule of $\coker(\alpha)$ is $T^\OO_q$;
\item[(b)] the torsion submodule of $\coker(\alpha^\ast)$ is canonically isomorphic to
$\e'(T^\OO_q)$;
\item[(c)] the vertical maps of 
diagram (\ref{(6-24)}) determine an isomorphism $T_q^\OO\to \e'(T^\OO_q)$, which
coincides with the linking form $\L^\OO$.
\end{itemize}
The star $\ast$ in diagram (\ref{(6-24)}) 
denotes the duality for finitely generated projective modules, cf. \cite{Fa3}, \S 1.4. For example $F_0^\ast$ denotes the set of all anti-linear homomorphisms $F_0\to \OO$. 
Note that in \cite{FL}, \S 1, we defined a cohomological version of the linking form $\L^\OO$;
in statement (c) above we refer to its homological analogue. 

Morphism $\alpha$ and $h$ in (\ref{(6-24)}) 
can be represented by matrices with entries in $\OO$. Hence,
we may find some $\epsilon>0$ such that all matrix elements of $\alpha$ and $h$ 
are analytic on $[-\epsilon, \epsilon]$. Hence, "tensoring" (\ref{(6-24)}) 
by $\ell^2(I_\epsilon, \mu)$, we
obtain the following commutative diagram
\begin{eqnarray}
\CD
(\ell^2(I_\epsilon, \mu))^{n_1}@>{\tilde \alpha}>> (\ell^2(I_\epsilon, \mu))^{n_0}\\
@V{(-1)^{q+1}\tilde h^\ast}VV @VV\tilde hV\\
(\ell^2(I_\epsilon, \mu))^{n_0}@>{\tilde \alpha^\ast}>>(\ell^2(I_\epsilon, \mu))^{n_1}
\endCD
\label{(6-25)}
\end{eqnarray}
where $n_i$ denotes the rank of $F_i$, $i=0, 1$. Here $\tilde \alpha$ and $\tilde h$ are given by the
same matrices as $\alpha$ and $h$, where instead of germs we have chosen analytic functions  $I_\epsilon\to \C$. 
More precisely, we use correspondence between functions $f:(-\epsilon, \epsilon)\to \C$
and functions $\tilde f: I_\epsilon\to \C$, where $\tilde f(c\exp(i\tau)) = f(\tau)$, 
$\tau\in [-\epsilon, \epsilon]$. 

We claim that for $\epsilon>0$ small enough:
\begin{itemize}
\item[(a')] the torsion submodule of $\coker(\tilde \alpha)$ is $T_q$, the torsion part of extended homology $\H_q(K;\ell^2(I_\epsilon, \mu))$;
\item[(b')] the torsion submodule of $\coker(\tilde \alpha^\ast)$ is canonically isomorphic to
$\e(T_q)$;
\item[(c')] the vertical maps of diagram (\ref{(6-25)}) determine an isomorphism $T_q\to \e(T_q)$, which
coincides with the linking form $\L$.
\end{itemize}
Here we view $\tilde \alpha$ and $\tilde \alpha^\ast$ as morphisms of the extended category, and
their cokernels are understood in this sense.

The statements (a'), (b'), (c') follow from the discussion of section \ref{sec3.7} 
(passage from diagram (\ref{(3-18)}) to diagram (\ref{(3-19)}).

Let us show that {\it for small enough $\epsilon$ there is a canonical isomorphism 
$F_c(T_q^\OO)\to T_q$ and via this isomorphism the form
$\L$ turns into the push forward form of $\L^\OO$}. 

To show this, first note that without loss of generality
we may assume that in diagram (\ref{(6-24)}) $\OO$-homomorphism $\alpha$ is injective.
Indeed, if it is not injective we may represent $F_1$ as a direct sum $K\oplus F_1'$, where
$K$ is the kernel of $\alpha$ and $\alpha|_{F_1'}$ is injective. Then $\alpha^\ast$ takes values in
$(F_1')^\ast$. Therefore, replacing $F_1$ by $F_1'$ we arrive at a similar diagram with injective $\alpha$. 

Let $F_0'\subset F_0$ be the smallest $\OO$-submodule containing $\im(\alpha)$ such that $F_0/F_0'$ 
has no torsion. Then $F_0\simeq F_0'\oplus (F_0/F_0')$ and $\alpha^\ast$ vanishes on $(F_0/F_0')^\ast$.
Hence the diagram 
\begin{eqnarray*}
\CD
F_1@>{\alpha}>> F_0'\\
@V{(-1)^{q+1}h^\ast}VV @VVhV\\
{F_0'}^\ast@>{\alpha^\ast}>> F_1^\ast
\endCD\end{eqnarray*}
will replace (\ref{(6-24)}) and here $\alpha$ is injective and $\coker(\alpha)$ is torsion. 

Both above changes (cutting off the kernel and the free part of the cokernel of $\alpha$) on diagram
(\ref{(6-24)}) can be performed simultaneously on the corresponding diagram (\ref{(6-25)}) assuming that $\epsilon$
is small enough. Both changes do not influence the resulting torsion forms. Hence we may assume
that the initial diagram (\ref{(6-24)}) has injective $\alpha$ and $\coker(\alpha)$ is torsion. But then our 
italicized above statement follows directly from the definitions. 

This completes the proof of Theorem \ref{claim}. \qed

\subsection{End of proof of Theorems \ref{thm6.6} and \ref{thm6.7}} 
Consider first the case of $q$ odd, when $\L$ and $\L^\OO$ are Hermitian. We will assume that
the interval $J$ contains $c$ as an internal point as assumed in Theorem \ref{thm6.6}.

Because of Theorem \ref{claim} and Remark 2.12 in \cite{FL}, the form $\L$ 
can be represented as an orthogonal sum of forms $\L_{k, \pm}$, which are defined as 
the discriminant form of 
\begin{eqnarray}
(\pm t^k: \ell^2([-\epsilon, \epsilon],\mu)\to \ell^2([-\epsilon, \epsilon],\mu)).
\label{(6-27)}
\end{eqnarray}
Here $\mu$ is the Lesbegue measure on the interval $[-\epsilon, \epsilon]$ and $\pm t^k$ stands
for the multiplication operator by the function $\pm t^k$. It is clear that
\begin{eqnarray}
\c_+(\L_{k,+}) = k,\quad \c_-(\L_{k,+})
=\cases
k, \quad\text{if $k$ is odd}\\
0, \quad\text{if $k$ is even}
\endcases\label{(6-28)}
\end{eqnarray} 
and similarly
\begin{eqnarray}
\c_-(\L_{k,-}) = k,\quad \c_+(\L_{k,-})
=\cases
k, \quad\text{if $k$ is odd}\\
0, \quad\text{if $k$ is even}
\endcases\label{(6-29)}
\end{eqnarray}
Statements (\ref{(6-28)}) and (\ref{(6-29)}) 
are based on the calculation similar to that made in Theorem 4.13 in \cite{Fa2}.
If the form $\L$ contains $n_k^\pm(c)$ copies of the form $\L_{k, \pm}$ we obtain that
$\c_+(\L)$ equals the maximum of numbers $k$, such that $n_k^+>0(c)$ or 
$k$ is odd and $n_k^->0(c)$ (by Proposition \ref{prop6.4}).

This proves Theorem \ref{thm6.6} in case $q$ odd; 
here we use Proposition 9.8 from \cite{FL}, expressing the linking form of deformation 
(defined in \S 1 of \cite{FL}, cf. formula (6)) in terms of the Blanchfield form. 

Assuming that the interval $J$ contains $c$ as its terminal point (as in Theorem \ref{thm6.7})
we will have instead of (\ref{(6-28)}), (\ref{(6-29)})
\begin{eqnarray*}
\c_+(\L_{k,+}) = \c_-(\L_{k,-})
=\cases
k, \quad\text{if $k$ is even}\\
0, \quad\text{if $k$ is odd}
\endcases
\end{eqnarray*}
and 
\begin{eqnarray*}
\c_+(\L_{k,-}) = \c_-(\L_{k,+})
=\cases
0, \quad\text{if $k$ is even}\\
k, \quad\text{if $k$ is odd}
\endcases
\end{eqnarray*}
Taking maximum and using Proposition \ref{prop6.4} leads to (\ref{(6-12)}).

The case when $q$ is even is similar. \qed

\section{\bf  The torsion signature}\label{sec7}

In this section we introduce another numerical invariant of torsion Hermitian forms - the torsion signature.

\subsection{} Let $\CC$ be a superfinite von Neumann category and let $\tr$ be a fixed
not normal (i.e. Dixmier type) 
trace on $\CC$. We refer to \cite{Fa2}, \cite{Fa4}, where these notions are described
in detail. 

In \cite{Fa4} we introduced the notion of {\it torsion dimension}
$$\tdim (\X)\in \R$$ 
for
any torsion object $\X\in \ob(\eca)$ with respect to the trace $\tr$. 
This concept of dimension for torsion objects behaves similarly 
to the well known von Neumann dimension for projective objects.
In particular, 
\begin{eqnarray}\tdim(\X\oplus \Y) = \tdim(\X) +\tdim(\Y).\end{eqnarray}
Compare the property
\begin{eqnarray}\c(\X\oplus \Y) =\max\{\c(\X), \c(\Y)\}\end{eqnarray}
for the Novikov-Shubin capacity. 

Recall that $\tdim (\X)$ is defined as the infimum of $\dim P$ (with respect to the trace $\tr$),
such that $P$ is projective and can be mapped onto $\X$. 
The torsion dimension $\tdim(\X)$
can also be characterized and computed in terms of the spectral density function, cf. \cite{Fa4}.

Using the notion of torsion dimension, we may define a new numerical invariants of torsion Hermitian
forms $\L: \X\to \e(\X)$. Namely, by Theorem 7.7 of \cite{Fa3} we can uniquely split 
$\X = \X_+ \oplus \X_-$ as the orthogonal sum of positively and
negatively definite forms. Now we set
\begin{eqnarray}
\tsig(\L) = \tsig_{\tr} (\L) = \tdim(\X_+) - \tdim(\X_-).
\label{(7-3)}
\end{eqnarray}
We will call number (\ref{(7-3)}) {\it the torsion signature} of form $\L$.

It is clear, that {\it the torsion signature of a hyperbolic form is zero}. 

It is certainly
not true that the torsion signature of a metabolic form vanishes 
(since, as we know from \cite{Fa3},\S 5,
all torsion forms are metabolic).

Note also that the torsion signature of an orthogonal sum $\L_1\perp \L_2$ equals the sum of their
torsion signatures
\begin{eqnarray}
\tsig(\L_1\perp \L_2) = \tsig(\L_1)+\tsig(\L_2).\label{(7-4)}
\end{eqnarray}
\subsection{Torsion signatures of odd-dimensional manifolds} Let $K$ be a closed
$(2q+1)$-dimensional closed oriented manifold and let $\E\to K$ be a unitary flat bundle
with fiber and monodromy in a superfinite von Neumann category $\CC$ supplied with a 
choice of a Dixmier type trace. 
By Theorem \ref{thm3.4}, there is a well-defined
$(-1)^{q+1}$-Hermitian non-degenerate linking form
\begin{eqnarray*}\L: T(\H_q(K;\E)) \to \e(T(\H_q(K;\E))).\end{eqnarray*}
We will define {\it the torsion signature of $K$} as follows
\begin{eqnarray*}\tsig_{\tr}(K;\E) = \cases \tsig_{\tr}(\L),\quad\text{if}\quad q\quad\text{is odd,}\\
\tsig_{\tr}(- \sqrt{-1}\cdot \L),\quad\text{if}\quad q\quad\text{is even.}
\endcases\end{eqnarray*}

The torsion signature is a homotopy invariant of $K$.

Changing the orientation of $K$ reverses the sign of the torsion signature:
\begin{eqnarray*}\tsig_{\tr}(-K;\E) = - \tsig_{\tr}(K;\E).\end{eqnarray*}

There is always a canonical choice for the flat bundle $\E$; namely, we may take for $\E$ the 
bundle with fiber $\ell^2(\pi_1(K))$ with the standard action of the fundamental group $\pi_1(K)$. 

Theorem 5.1 gives a necessary condition for vanishing of the torsion signature: it happens if our
odd-dimensional manifold can be obtained as a slice of an even-dimensional manifold with
vanishing middle-dimensional extended homology. In section \ref{sec8} we will give a different
criterion for vanishing of the torsion signature.

\subsection{The torsion signature and the Blanchfield form} Here we will compute explicitly
the torsion
signature for closed oriented manifolds $K$, with $\dim K = 2q+1$, $\pi_1(K) =\Z$; 
we will take for the flat bundle $\E$ the
natural flat bundle 
with fiber $\ell^2(\Z)=\ell^(S^1)$, cf. section \ref{sec6.5}. 

As in section \ref{sec6.5},
our von Neumann category
$\CC$ is the category of finite dimensional Hilbert spaces over the circle $S^1$. 
In Theorems \ref{thm6.6} and \ref{thm6.7} 
we considered two different traces on $\CC$ and both these traces were
normal. Now we will consider non-normal traces on $\CC$, which we will construct as follows.
Let $c\in \mathcal S_q(K)$ be a point of the unit circle, which belongs to the support of the torsion submodule of $H_q(\tilde K;\C)$, cf. \ref{sec6.5}. 
Let $\xi\mapsto \H(\xi)$ be a measurable field of finite
dimensional Hilbert spaces over $S^1$ and let $f: \H \to \H$ be a decomposable linear map, cf. \cite{Di}, part II, chapter 2. 
For almost all $\xi\in S^1$ we have a linear map $f(\xi):\H(\xi)\to \H(\xi)$ and we will
denote by $\text{Tr}(f(\xi))$ its usual trace (the sum of diagonal entries in a matrix representation).
We will define now two traces $\tr_\pm$ on category $\CC$ by setting
\begin{eqnarray}
\tr_+(f) =\lo \, \, [n\int_0^{1/n}\text{Tr}(f(c\exp(i\phi)))d\phi]\label{(7-8)}
\end{eqnarray}
and   
\begin{eqnarray}\tr_-(f) =\lo \, \, [n\int_{-1/n}^0\text{Tr}(f(c\exp(i\phi)))d\phi]
\label{(7-9)}
\end{eqnarray}   
Note that the function $\xi\mapsto \text{Tr}(f(\xi)))$, where $\xi\in S^1$,
is essentially bounded; hence 
$$n\mapsto [n\int_0^{1/n}\text{Tr}(f(c\exp(i\phi))d\phi]$$ 
is a bounded sequence on complex numbers, and therefore we may apply in (\ref{(7-8)}) 
the Banach limit $\lo$ as $n\to\infty$. 
The same explanation applies to (\ref{(7-9)}). Formulae (\ref{(7-8)}) and (\ref{(7-9)}) define
two different non-normal traces on $\CC$. Both these traces are "supported" at point $c\in S^1$, 
but "they live on different sides of $c$". 

\begin{theorem}\label{thm8.1} Let $K$ be a closed oriented $(2q+1)$-dimensional manifold with
$\pi_1(K)=\Z$, and let
$\L: T_q\to \e(T_q)$, where $T_q=T(\H_q(K;\ell^2(\Z)))$, be the linking form. Then the torsion
signature of $\L$ with respect to the trace $\tr_+$ equals
\begin{eqnarray}\tsig_{\tr_+}(\L) \, =\, \sigma_{\text{ev}}(c) +\sigma_{\text{odd}}(c),\end{eqnarray}
and the torsion
signature of $\L$ with respect to the trace $\tr_-$ equals
\begin{eqnarray}\tsig_{\tr_-}(\L)  \, = \,
\sigma_{\text{ev}}(c) -\sigma_{\text{odd}}(c).\end{eqnarray}
Here 
\begin{eqnarray}\sigma_{\text{ev}}(c) = \sum_{k\, \text{even}} [n_k^+(c) - n_k^-(c)]\quad\text{and}\quad
\sigma_{\text{odd}}(c) = \sum_{k\, \text{odd}} [n_k^+(c) - n_k^-(c)],\end{eqnarray}
where $n_k^\pm(c)$ are the numbers determined by the Blanchfield form, cf. section \ref{sec6.8}.
\end{theorem}
\begin{proof} The proof is similar to the proof of Theorems \ref{thm6.6} and \ref{thm6.7} 
described above in full detail. It is based on Theorem \ref{claim}, 
which explicitly computes the linking form $\L$. 
The distinction happens only at the very last stage of the proof, when
one computes the torsion signatures of the elementary forms (\ref{(6-27)}) 
and uses additivity (\ref{(7-4)}). 
\end{proof}

Comparing Theorems \ref{thm6.6}, \ref{thm6.7} and \ref{thm8.1} 
we conclude, that {\it the torsion signature is 
more informative than the Novikov - Shubin signature.}

The numbers $\sigma_{\text{odd}}(c)$ and $\sigma_{\text{ev}}(c)$ play an important role
in \cite{FL}, where they describe the behavior of the eta-invariant.

\section{\bf Excess of extension}\label{sec8}

We know that the torsion signature of a hyperbolic form vanishes. In this section we will give a 
different criterion for vanishing of the torsion signature.

Let $\CC$ be a von Neumann category with a non-normal trace $\tr$. 
Everywhere in this section we will assume that the trace $\tr$ is non-negative, cf. \cite{Fa2}.

For any short exact sequence
\begin{eqnarray}0\to \X' \to \X \to \X'' \to 0,\label{(8-1)}\end{eqnarray}
consisting of torsion objects of the extended abelian category $\eca$, holds
\begin{eqnarray}\max \{ \tdim(\X'), \tdim(\X'')\}  \le \tdim (\X) \le \tdim (\X') + \tdim (\X''),\label{(8-2)}
\end{eqnarray}
cf. \cite{Fa4}. 
Moreover, if sequence (\ref{(8-1)}) splits, then
$\tdim (\X) = \tdim (\X') + \tdim (\X'').$
We may use inequality (\ref{(8-2)}) in order to measure
the complexity of extension (\ref{(8-1)}).

\begin{definition}
Given extension (\ref{(8-1)}), define {\it its excess} 
(with respect to the chosen trace) as the following
real number
\begin{eqnarray}
\ex(\X'\to \X\to \X'') = \tdim (\X') + \tdim (\X'') - \tdim (\X). \label{(8-3)}\end{eqnarray}
\end{definition}

Non triviality of the excess (with respect to at least one trace $\tr$ on $\CC$) implies, that
the extension does not split. Note that
\begin{eqnarray}0 \le \ex(\X'\to \X\to \X'') \le \min\{\tdim (\X'), \tdim (\X'')\}.
\label{(8-4)}\end{eqnarray}
Lemma \ref{l8.5} below provides many examples showing that the excess may be nonzero.

\begin{theorem}\label{thm8.2} If a non-degenerate torsion Hermitian form $\L: \X \to \e(\X)$
in $\eca$ admits a metabolizer $\Y \subset \X$, $\Y^{\perp} = \Y$, with vanishing excess
\begin{eqnarray}\ex(\Y \to \X \to \X/\Y) = 0,\label{(8-5)}\end{eqnarray}
then the torsion signature of $\L$ vanishes:
\begin{eqnarray}\tsig(\L) = 0.\label{(8-6)}\end{eqnarray}
\end{theorem}

The proof of Theorem \ref{thm8.2} will be given later at the end of this section. 
The proof will use the following 
Lemmas.

\begin{lemma}\label{l8.3} If $\Y\subset \X$ is a metabolizer, then
vanishing of excess (\ref{(8-5)}) is equivalent to 
\begin{eqnarray}2\tdim(\Y) = \tdim(\X).\label{(8-7)}\end{eqnarray}
\end{lemma}
\begin{proof}
For any subobject $\Y \subset \X$ we have
\begin{eqnarray}\tdim(\X/\Y) = \tdim(\Y^{\perp})\label{(8-8)}\end{eqnarray}
because of the canonical isomorphism $\Y^{\perp} \simeq \e(\X/\Y)$, and $\e(\X/\Y)$ is
isomorphic to $\X/\Y$ (not canonically). If $\Y$ is the metabolizer, then $\tdim(\X/\Y) = \tdim(\Y)$, (from (\ref{(8-8)})) and hence vanishing of excess (\ref{(8-5)}) is equivalent to (\ref{(8-7)}). 
\end{proof}

\begin{lemma}\label{l8.5} Suppose that $\L: \X \to \e(\X)$ is 
a non-degenerate positive definite torsion Hermitian form. 
Then for the metabolizer $\Y \subset \X$ (which is unique by
Proposition 6.2 of \cite{Fa3}) holds
\begin{eqnarray}\ex(\Y \to \X \to \X/\Y) = \tdim \X = \tdim \Y = \dim \X/\Y,\label{(8-9)}
\end{eqnarray}
i.e. this extension has the maximal possible excess. 

Any subobject $\mathcal Z \subset \X$, such that
$\mathcal Z^{\perp} \subset \mathcal Z$, has torsion dimension equal to torsion dimension of $\X$,
$\tdim \mathcal Z =\tdim \X$.
\end{lemma}
\begin{proof} From the proof of Proposition 6.2 in \cite{Fa3} we know that without
loss of generality we may assume that
$\X$ is represented as $(\alpha: A \to A)$, where $A$ is an object of $\CC$ and $\alpha$ is
a self-adjoint positive operator, $\alpha\in \hom_{\CC}(A, A)$, and the metabolizer
$\Y\subset \X$ is represented as $(\beta:A \to A)$ with $\beta$ being the positive
square root of $\alpha$. If $F(\lam)$ is the spectral density function of $\X$
then $F({\lam}^2)$ is the spectral density function of $\Y$. Therefore we obtain
$$\tdim \X = \lim_{\lam \to +0} F(\lam) = \lim_{\lam \to +0} F({\lam}^2) = \tdim \Y.$$

Now, formula 
(\ref{(8-8)}) shows that 
$$\ex(\Y \to \X \to \X/\Y) \, = \, 2\tdim(\Y)-\tdim(\X),$$
 which proves 
(\ref{(8-7)}).

If $\mathcal Z^{\perp} \subset \mathcal Z$ then $\mathcal Z$ contains the metabolizer of $\X$
(by the arguments in the proof of Propositions 6.1 and 6.2 in \cite{Fa3}) and
therefore $\tdim(\mathcal Z) = \tdim(\X)$. 
This completes the proof.  
\end{proof}

The following Corollary is a partial result in the direction of Theorem \ref{thm8.2}; it takes care
of the case of reduction of a positively defined form with respect to a submodule
with zero excess.

\begin{corollary} Let $\L: \X \to \e(\X)$ be a non-degenerate positive definite
torsion Hermitian form. Suppose that $\Y\subset \X$ is a subobject such that
$\Y \subset \Y^{\perp}$. Then the following conditions are equivalent:
\begin{itemize}
\item[(i)] the excess $\ex(\Y \to \X\to \X/\Y) = 0$ vanishes;
\item[(ii)]" 
$\tdim(\Y) = 0.$
\end{itemize}
Either of these conditions implies
\begin{eqnarray}\tdim(\X) = \tdim(\Y^{\perp}/\Y).\label{(8-10)}\end{eqnarray}
\end{corollary}
\begin{proof} If $\Y \subset \Y^{\perp}$ then by  Lemma \ref{l8.5}
$$\tdim(\Y^{\perp}) = \tdim(\X).$$
Thus, if the excess (i) vanishes, we get
$$
\begin{array}{l}
\tdim(\X) = \tdim(\Y) + \tdim(\X/\Y) \\
= \tdim(\Y) + \tdim(\Y^{\perp})\\
= \tdim(\Y) + \tdim(\X),
\end{array}
$$
which implies $\tdim(\Y) = 0$.

Conversely, if $\tdim(\Y) = 0$, then inequality (\ref{(8-4)}) shows that
$\ex(\Y \to\X \to \X/\Y) =0$.

If (ii) holds, then $\ex(\Y\to\Y^{\perp}\to\Y^{\perp}/\Y) =0$ (again, because of (\ref{(8-4)}))
and
$$\tdim(\Y^{\perp}/\Y) = \dim(\Y^{\perp}) = \tdim(\X).$$
This completes the proof. 
\end{proof}

\subsection*{Proof of Theorem \ref{thm8.2}} Let $\X= \X_+ \oplus \X_-$ be the canonical decomposition
determined by $\L$ (cf. \cite{Fa3}, Theorem 7.7). 
Denote by $\Y_\pm = i_\pm(\Y)\subset \X_\pm$ the image of $\Y$ under the projections
$j_\pm: \X\to \X_\pm$.
Then clearly $\Y^{\perp}_\pm \subset \Y_\pm$ (otherwise we will get a contradiction
to $\Y^{\perp} = \Y$).
Applying the second statement of Lemma \ref{l8.5} we find
\begin{eqnarray}\tdim(\Y_\pm) = \tdim(\X_\pm).\label{(8-11)}\end{eqnarray}
On the other hand,
\begin{eqnarray}\tdim(\Y_\pm) \le \tdim(\Y).\label{(8-12)}\end{eqnarray}
Therefore we obtain
\begin{eqnarray}\tdim(\X_+) \le \tdim(\Y),\label{(8-13)}\end{eqnarray}
\begin{eqnarray}\tdim(\X_-) \le \tdim(\Y).\label{(8-14)}\end{eqnarray}
Now we will use $\ex(\Y\to\X\to \X/\Y)=0$ to conclude
$$
2\tdim(\Y) = \tdim(\X) = \tdim(\X_+) + \tdim(\X_-)
$$
This shows that inequalities (\ref{(8-13)}) and (\ref{(8-14)}) are in fact equalities,
and therefore
$$\tdim(\X_+) = \tdim(\X_-)\quad \text{and}\quad \tsig(\L) =0.$$
This completes the proof. \qed

\bibliographystyle{amsalpha}

\end{document}